\documentclass[12pt,a4paper]{amsart}
\usepackage[applemac]{inputenc}
\usepackage[T1]{fontenc}
\usepackage{amssymb}
\usepackage{amsmath}
\usepackage{amsthm}
\usepackage{varioref} % pour l'enumerate avec lettres incorporés
\usepackage[pdftex]{hyperref}  % als laatste toevoegen

\usepackage{color}

\renewcommand{\theenumi}{\roman{enumi}}

\newtheoremstyle{mythmstyle}%
{}% space above
{}% space below
{\itshape}% bodyfont
{}% indentamount
{\bf}% theorem head font
{}% punctuation after thm head
{0pt}% space after thm head
{} % theorem head spec

\newtheoremstyle{mydefstyle}%
{}% space above
{}% space below
{}% bodyfont
{}% indentamount
{\bf}% theorem head font
{}% punctuation after thm head
{0pt}% space after thm head
{} % theorem head spec

\newtheoremstyle{mypreuvestyle}%
{}% space above
{}% space below
{}% bodyfont
{}% indentamount
{\em}% theorem head font
{}% punctuation after thm head
{0pt}% space after thm head
{} % theorem head spec

\newcounter{remembersection}
\newif\ifmafirstsection\mafirstsectiontrue

\newif\ifmynonumberenvi\mynonumberenvitrue

\theoremstyle{mythmstyle}
\newtheorem{proclaimmythm}{}
\numberwithin{proclaimmythm}{section}
\newtheorem*{proclaimmythm*}{}
\newenvironment{proclaim}[2][*]{\ifx*#1\mynonumberenvitrue\begin{proclaimmythm*}{\bf#2.} \ignorespaces\else\mynonumberenvifalse\begin{proclaimmythm}{. \bf#2.}\label{#1} \ignorespaces\fi}{\ifmynonumberenvi\end{proclaimmythm*}\else\end{proclaimmythm}\fi}

\theoremstyle{mydefstyle}
\newtheorem{proclaimmydef}[proclaimmythm]{}
\newtheorem*{proclaimmydef*}{}
\newenvironment{definition}[2][*]{\ifx*#1\mynonumberenvitrue\begin{proclaimmydef*}{\bf#2.}\else\mynonumberenvifalse\begin{proclaimmydef}{. \bf#2.}\label{#1} \ignorespaces\fi{}}{\ifmynonumberenvi\end{proclaimmydef*}\else\end{proclaimmydef}\fi}

\newif\ifmynonumberequation\mynonumberequationtrue

\makeatletter
\numberwithin{equation}{section}
\newenvironment{moneq}[1][*]{\ifx*#1\mynonumberequationtrue\begin{displaymath}\else\mynonumberequationfalse\begin{equation}\label{#1}\fi}{\ifmynonumberequation\end{displaymath}\@ignoretrue\else\end{equation}\@ignoretrue\fi\ignorespaces}
\makeatother

\def\QEDbox{\hbox{\lower2.3pt\vbox{\hrule\hbox
   {\vrule\kern1pt\vbox{\kern1.7pt\hbox{$\scriptstyle
   QED$}\kern.6pt}\kern1pt\vrule}\hrule}}}
\def\QED{\hskip0.01em plus 40pt\null{} \null\nobreak\hfill
   \kern3pt\QEDbox}

\theoremstyle{mypreuvestyle}
\newtheorem*{proclaimmypreuve}{}
\newenvironment{preuve}[1][{Proof}]{\begin{proclaimmypreuve}{\em#1.} \ignorespaces}{\QED\end{proclaimmypreuve}}

\newif\ifdiscussion\discussionfalse

\newcommand\body{{\mathbf B}}
\newcommand\CA{\mathcal{A}}
\newcommand\conjug{\mathfrak{C}}
\newcommand\contrf[2]{\iota(#1)#2}
\newcommand\End{{\operatorname{End}}}
\newcommand\extder{{\operatorname{d}}}
\newcommand\fracp[2]{\frac{\partial#1}{\partial#2}}
\newcommand\gammab{{\bar\gamma}}
\newcommand\Gammah{{\widehat\Gamma}}
\newcommand\Gammat{{\widetilde\Gamma}}
\newcommand\gammat{{\widetilde\gamma}}

\newcommand\hidecomputations[1]{}
\newcommand\ie{i.e.}
\newcommand\nablah{{\widehat\nabla}}
\newcommand\pb{{\bar p}}
\newcommand\Phih{{\widehat\Phi}}
\newcommand\Phit{{\widetilde\Phi}}
\newcommand\pr{{\operatorname{pr}}}
\newcommand\Psih{{\widehat\Psi}}
\newcommand\Psit{{\widetilde\Psi}}
\newcommand\recallp[1]{[\ref{#1}]}
\newcommand\recalle[1]{(\ref{#1})}
\newcommand\RR{\mathbf{R}}
\newcommand\scirc{\,{\raise 0.8pt\hbox{$\scriptstyle\circ$}}\,}
\newcommand\tanvec{{\mathcal V}}
\newcommand\vb{{\bar v}}

\allowdisplaybreaks%[3]

\begin{document}

% PERSONAL DATA
\author{Thomas Leuther}
\address[Thomas Leuther]{University of Li\`ege, Department of Mathematics\\
 Grande Traverse, 12 - B37, B-4000 Li\`ege, Belgium}
\email{Thomas.Leuther[at]ulg.ac.be}
\author{Fabian Radoux} 
\address[Fabian Radoux]{University of Li\`ege, Department of mathematics \\
 Grande Traverse, 12 - B37, B-4000 Li\`ege, Belgium}
\email{Fabian.Radoux[at]ulg.ac.be}
\author{Gijs M. Tuynman}
\address[Gijs M. Tuynman]{Laboratoire Paul Painlev\'e, U.M.R. CNRS 8524 et UFR de Math\'ematiques,
Universit\'e de Lille I, 59655 Villeneuve d'Ascq Cedex, France}
\email{Gijs.Tuynman[at]univ-lille1.fr}

% THIS DOCUMENT
%\date{\today} 

\title[Geodesics on a supermanifold]{
Geodesics on a supermanifold and projective equivalence of super connections
%Projective equivalence of super connections and their geodesics
}

\begin{abstract}
We investigate the concept of projective equivalence of connections in supergeometry. To this aim, we propose a %new definition for (super) geodesics on a supermanifold. As in the classical case, our (super) geodesics 
definition for (super) geodesics on a supermanifold in which, as in the classical case, they are the projections of the integral curves of a vector field on the tangent bundle: the \emph{geodesic vector field} associated with the connection. 
Our (super) geodesics possess the same properties as the in the classical case: there exists a unique (super) geodesic satisfying a given initial condition and when the connection is metric, our supergeodesics coincide with the trajectories of a free particle with unit mass. 
Moreover, using our definition, we are able to establish Weyl's characterization of projective equivalence in the super context: two torsion-free (super) connections define the same geodesics (up to reparametrizations) if and only if their difference tensor can be expressed by means of a (smooth, even, super) $1$-form.
\end{abstract}

%%%
\maketitle



{\bf Subject classification:} supermanifolds and supergroups, real and complex differential geometry.
%%%

\maketitle

%\tableofcontents

%%%%%%%%%%%%%%%%%%%%%%%%%%%%%%%%%%%%%%%%% <INTRODUCTION>
\section{Introduction}

%\textcolor{blue}{Il y a un problème avec le titre dans le sens que c'est trop long pour mettre au-dessus les pages (le numéro de page se superpose en 11pt et en 12pt ça dépasse largement). On pourrait laisser la décision à la revue où on le soumet, mais on peut aussi prendre le devant. On pourrait changer en "Projective equivalence of super connections and their geodesics". Je n'ai pas d'opinion précise sur ce sujet.}%

% caractÈrisation de Weyl classique :

%% - historique

The concept of projective equivalence of connections goes back to the 1920's, with the study of the so-called ``geometry of paths'' (see \cite{Th, TV, Wh} or \cite{Ro1,Ro2,HR} for a modern formulation). In 2002, M.~Bordemann used this theory to answer the problem of projectively invariant quantization in \cite{Bo}.

%% - lien avec les quantifications Èquivariantes

Projectively invariant quantization is a generalization to arbitrary manifolds of the notion of equivariant quantizations in the sense of Lecomte-Ovsienko, see \cite{LO,L,MR}. It consists in building in a natural way a quantization (i.e., a symbol-preserving linear bijection between a space of symbols and a space of differential operators) from a linear connection, requiring that the quantization remains unchanged if we start from another connection in the same projective class.

%% - dÈfinition gÈomÈtrique de l'Èquivalence projective et caractÈrisation algÈbrique

By definition, two connections are called {\em projectively equivalent} if they have the same geodesics, up to parametrization. In other words, the geodesics of two equivalent connections are the same, provided that we see them as sets of points, rather than as maps from an open interval of $\mathbf{R}$ into the manifold. In \cite{We}, H. Weyl showed that projective equivalence can be rephrased in an algebraic way: two connections are projectively equivalent if and only if the symmetric tensor which measures the difference between them can be expressed by means of a $1$-form.

% problËme de la supÈrisation :

%% - supÈrisation via la caractÈrisation algÈbrique + super QPI

Weyl's algebraic characterization of projective equivalence provides a convenient way to transport projective equivalence to the framework of supergeometry: two superconnections are said to be {\em projectively equivalent} if the (super)symmetric tensor which measures the difference between them can be expressed by means of a (super)$1$-form. Using this notion, it is possible to set the problem of projectively invariant quantization on supermanifolds while M.~Bordemann's method can be adapted in order to solve it (see \cite{LR}).

Remembering the classical picture, it is natural to ask whether it is possible to find a geometric counterpart to the algebraic definition of projective equivalence of superconnections, i.e., a characterization in terms of supergeodesics. The main purpose of the present paper is to answer this question in the affirmative.

%% - dÈfinition des gÈodÈsiques

As in the classical case, we define, in section \ref{sectionSuperGeodesics}, supergeodesics associated with a superconnection $\nabla$ on a supermanifold $M$ as being the projections onto $M$ of the integral curves of a vector field $G$ on the tangent bundle $TM$: the {\em geodesic vector field} of $\nabla$. In section \ref{proj} we then define the notion of reparametrization of a geodesic and establish that two connections $\nabla$ and $\nablah$ on a supermanifold $M$ have the same geodesics up to parametrization if and only if there is an even $1$-form $\alpha$ such that
\[
\nablah_{X}Y=\nabla_{X}Y+\alpha(X)Y+(-1)^{\epsilon(X)\epsilon(Y)}\alpha(Y)X\quad\forall X,Y\in\Gamma(TM) ,
\] 
thus showing that Weyl's characterization also holds in supergeometry.

% mise en perspective de notre approche de la notion de gÈodÈsique

%% - "a misconception concerning super curves"

%% vs Goertsches

We note that our approach to supergeodesics differs from that of Goertsches \cite{Go}. In particular, our equations for supergeodesics are the natural generalization of the classical ones. Actually, our approach is nearly identical to that recently proposed by Garnier-Wurzbacher in \cite{GW}, where they consider supergeodesics associated with a Levi-Civita superconnection. In their paper, supergeodesics on a Riemannian supermanifold $M$ are shown to co{\"\i}ncide with the projections of the flow of a Hamiltonian supervector field defined on the (even) cotangent bundle of $M$. In section \ref{Smetricconnections} we will show that the same holds in our approach when we use a Levi-Civita connection. 
In fact, beyond the fact that they restrict to the Riemannian setting, the main difference between Garnier-Wurzbacher's supergeodesics and ours lies in the domain of supercurves. Supercurves should be images of $1$-dimensional manifolds, but as it is well-known, the theory of supercurves with a single parameter turns out to be very shallow: supercurves in a single even parameter are reduced to ordinary curves in the body of the manifold while supercurves in a single odd parameter are simply odd straight lines. In order to overcome these limitations, we choose to change the viewpoint. Usually curves do not come singly, they appear in families. And in particular the integral curves of a vector field on a supermanifold $N$ should not be seen as a simplistic collection of curves, but as a map (the flow) defined on (an open subset of) $\RR \times N\,$(\footnote{{In fact, rather $\CA_0 \times N$ than $\RR \times N$ since maps defined on $\CA_0 \times N$ live in the category of supermanifolds while containing the same information as maps defined on $\RR \times N$ (see Lemma \ref{expandinginoddpowers}).}}), incorporating the initial condition in the domain of the map. And indeed, the flow of a vector field is jointly smooth in the time parameter $t$ and the initial condition $n\in N$. In the simplistic viewpoint one writes $\gamma_{n}(t)$ for an integral curve with initial condition $n\in N$, whereas in the viewpoint of a flow one rather writes $\varphi_t(n)$ or $\varphi(t,n)$. Roughly speaking, we could say that our change of viewpoint enlarges in a natural way (we do not add an arbitrary manifold $S$ as in \cite{GW}) the domain of supercurves so that it is now possible to get supercurves with desirable properties.

%%%%%%%%%%%%%%%%%%%%%%%%%%%%%%%%%%%%%%%%% </INTRODUCTION>

\section{Notation and general remarks}

We will work with the geometric $H^\infty$ version of DeWitt supermanifolds, which is equivalent to the theory of graded manifolds of Leites and Kostant (see \cite{DW, Ko, Le, Rog, Tu1}). Any reader using a (slightly) different version of supermanifolds should be able to translate the results to her\slash his version of supermanifolds.

\begin{definition}{Some general conventions}

\begin{itemize}
\item
The basic graded ring will
be denoted as $\CA$ and we will think of it as the exterior algebra $\CA
= \Lambda V$ of an infinite dimensional real vector space $V$.

\item
Any element $x$ in a graded space splits into an even and an odd part $x=x_0 + x_1$. {Associated to this splitting we have the operation $\conjug$ of conjugation in the odd part defined by $\conjug(x) \equiv \conjug(x_0+x_1) = x_0 - x_1$.} %The parity function $\e$ is defined on homogeneous elements, \ie, elements $x$ for which either the even part $x_0$ or the odd part $x_1$ is zero. More precisely, if $x=x_\alpha$,  then $\e(x) = \alpha$.

\item
All (graded) objects over the basic ring $\CA$ have an underlying real
structure, called their \emph{body}, in which all nilpotent elements in $\CA$ are
ignored\slash killed. This forgetful map is called the \emph{body map}, denoted by $\body$. 
{For the ring $\CA$, this map $\body$ is nothing but the canonical projection $\CA =\Lambda V \to \Lambda^0 V = \RR$.}

\item
If $\omega$ is a $k$-form and $X$ a vector field, we denote the contraction of
the vector field $X$ with the $k$-form $\omega$ by $\contrf{X}{\omega}$, which
yields a $k-1$-form. If $X_1, \dots, X_\ell$ are $\ell\le k$ vector fields, we
denote the repeated contraction of $\omega$ by $\contrf{X_1, \cdots,
X_\ell}{\omega}$. More precisely:
\begin{moneq}
\contrf{X_1, \cdots, X_\ell}{\omega} = 
\Bigl( \contrf{X_1}{} \scirc \cdots \scirc \contrf{X_\ell}{}
\Bigr) \, \omega
\end{moneq}%
In the special case $\ell = k$ this definition differs by a factor
$(-1)^{k(k-1)/2}$ from the usual definition of the evaluation of a $k$-form on
$k$ vector fields. This difference is due to the fact that in ordinary
differential geometry repeated contraction with $k$ vector fields corresponds
to the direct evaluation in the reverse order. And indeed, $(-1)^{k(k-1)/2}$ is
the signature of the permutation changing $1,2,\dots,k$ in $k,k-1,\dots,2,1$.
However, in graded differential geometry this permutation not only introduces
this signature, but also signs depending upon the parities of the vector
fields. These additional signs are avoided by our definition.

\item
Evaluation\slash contraction of a left-(multi-)linear map $f$ with a vector $v$ is denoted just as the contraction of a differential form with a vector field as $\contrf{v}{f}$. If $f:E\to \CA$ is just left-linear, this is just the image of $v$ under the map $f$. However, if $f$ is for instance left-bilinear, the contraction $\contrf{v}{f}$ now is a left-linear map given by
\begin{moneq}
\contrf v f : w\mapsto \contrf{w,v}f
\end{moneq}%
As left-linearity and right-linearity are the same for even maps, we sometimes use the more standard notation $f(w,v)$ for the image of the couple $(w,v)$ under the bilinear map $f$, instead of $\contrf{w,v}f$.

\item
If $E$ is an $\CA$-vector space, $E^*$ will denote the {\em left} dual of $E$, \ie, the
space of all {\em left}-linear maps from $E$ to $\CA$.

\item
Let $x^1, \dots, x^n$ be local coordinates of a super manifold $M$ of graded dimension $p\vert q$, $p+q=n$, ordered such that $x^1, \dots, x^p$ are even and $x^{p+1}, \dots, x^n$ are odd (we will denote the latter also by $(\xi^1, \dots, \xi^q)$). Using the symbol $\varepsilon$ as the parity function, we thus have $\varepsilon(x^i)=0$ for $i\le p$ and $1$ for $i>p$. To simplify notation, we introduce the abbreviation $\varepsilon_i = \varepsilon(x^i)$.

\end{itemize}

\end{definition}

\begin{proclaim}[expandinginoddpowers]{Lemma (\cite{Tu1})}
Let $f$ and $g$ be smooth functions of even variables $x_1, \dots, x_{{p}}$ and odd variables $\xi_1, \dots, \xi_{{q_1}}$ and $\eta_1, \dots, \eta_{{q_2}}$. We can expand these functions with respect to products of odd variables, either only the $\xi$'s, only the $\eta$'s or both $\xi$'s and $\eta$'s, giving (for $f$) the formulae
\begin{align*}
f(x,\xi,\eta) &= \sum_{I\subset \{1, \dots, {q_1}\}} \xi^I \cdot f_I^{(\xi)}(x,\eta)
=
\sum_{J\subset\{1, \dots, {q_2}\}} \eta^J \cdot f_J^{(\eta)}(x,\xi)
\\&
=
\sum_{I\subset \{1, \dots, {q_1}\}, J\subset\{1, \dots, {q_2}\}} \xi^I \cdot \eta^J \cdot f_{IJ}^{(\xi,\eta)}(x)
\end{align*}
where the sum is over all subsets with (for instance)
\begin{moneq}
I= \{i_1, \dots, i_k\} \text{ with } 1\le i_1 < i_2 < \cdots <i_k\le {q_1} 
\qquad\Longrightarrow\qquad
\xi^I = \xi_{i_1} \cdots \xi_{i_k}
\end{moneq}%
Then the following statements are equivalent:
\begin{enumerate}
\item $f=g$

\item for all $I\subset \{1, \dots, {q_1}\}$: $f_I^{(\xi)} = g_I^{(\xi)}$

\item for all $J\subset \{1, \dots, {q_2}\}$: $f_J^{(\eta)} = g_J^{(\eta)}$

\item for all $I\subset \{1, \dots, {q_1}\}$, $J\subset \{1, \dots, {q_2}\}$: $f_{IJ}^{(\xi,\eta)} = g_{IJ}^{(\xi,\eta)}$

\end{enumerate}
Moreover, when we have expanded with respect to all odd variables, the remaining functions of the even variables only are completely determined by their values on real coordinates. Said differently, we may assume that they are ordinary smooth functions of $n$ real coordinates.

\end{proclaim}

\section{Super Geodesics}
\label{sectionSuperGeodesics}

{Before dealing with the specific problem of geodesics on a supermanifold, we first recall some general definitions and facts about (super) connections in the tangent bundle. Then we attack the problem of defining super geodesics: we associate with any connection a so-called \emph{geodesic vector field} on the tangent bundle, whose flow equations are the straightforward super analogs of the classical geodesic equations.
}

\begin{definition}{Definition \cite[VII\S6]{Tu1}}
A connection in a (super) vector bundle $p:E\to M$ over a supermanifold $M$ is (can be seen as) a map $\nabla:\Gamma(TM) \times \Gamma(E) \to \Gamma(E)$ satisfying
\begin{enumerate}
\item $\nabla$ is bi-additive (in $\Gamma(TM)$ and $\Gamma(E)$) 

\item
for $X\in \Gamma(TM)$, $s\in \Gamma(E)$ and $f\in C^\infty(M)$ we have
\begin{moneq}
\nabla_{fX}s = f\cdot \nabla_Xs
\end{moneq}%

\item
for homogeneous $X\in \Gamma(TM)$, $s\in \Gamma(E)$ and $f\in C^\infty(M)$ we have
\begin{moneq}
\nabla_{X}(fs) = (Xf)\cdot s + (-1)^{\varepsilon(X)\cdot \varepsilon(f)} f\cdot \nabla_Xs
\end{moneq}%

\end{enumerate}

\end{definition}

\begin{proclaim}{Lemma}
If $\nabla$ and $\nablah$ are connections in $E$, the map $S:\Gamma(TM) \times \Gamma(E) \to \Gamma(E)$ defined by
\begin{moneq}
S(X,s) = \nabla_Xs - \nablah_Xs
\end{moneq}%
is even and bilinear over $C^\infty(M)${. In other words, $S$ is} a ``tensor'', \ie, can be seen as a section of the bundle $TM^* \otimes \End(E)$ \cite[IV\S5]{Tu1}.

\end{proclaim}

\begin{proclaim}{Lemma}
If $\nabla$ is a connection in $TM$, then the map $T:\Gamma(TM) \times \Gamma(TM) \to \Gamma(TM)$ defined on homogeneous $X,Y\in \Gamma(TM)$ by
\begin{moneq}
T(X,Y) = \nabla_XY - (-1)^{\varepsilon(X)\cdot \varepsilon(Y)}\cdot \nabla_YX - [X,Y]
\end{moneq}%
is even, graded anti-symmetric and bilinear over $C^\infty(M)${. In other words, $T$ is} a ``tensor'', \ie, can be seen as a section of the bundle $\bigwedge^2TM^* \otimes TM$, \ie, as a 2-form on $M$ with values in $TM$ \cite[IV\S5]{Tu1}.

\end{proclaim}

\begin{definition}{Definition}
A connection $\nabla$ in $TM$ is said to be {\emph{torsion-free}} if the tensor $T$ is identically zero.

\end{definition}

\begin{proclaim}{Corollary}
If $\nabla$ and $\nablah$ are {torsion-free} connections in $TM$, the tensor {$S = \nabla - \nablah: \Gamma(TM)\times \Gamma(TM)\to \Gamma(TM)$} is graded symmetric.

\end{proclaim}

 {Let $\nabla$ be a connection in $TM$ (we also say a connection \emph{on} $M$). {On a local chart for $M$ with coordinates $x=(x^1, \dots, x^n)$ we define the Christoffel symbols $\Gamma_{jk}^i$ of $\nabla$ by
\begin{moneq}
\Gamma_{jk}^i(x) 
=
\contrf{\nabla_{\partial_{x^j}}\partial_{x^k}}{\,\extder x^i\vert_x}
\end{moneq}%
%In particular, we have $\sum_i \Gamma_{jk}^i \cdot \partial_{x^i} = \nabla_{\partial_{x^j}}\partial_{x^k}$.
with parity $\varepsilon\bigl(\Gamma_{jk}^i(x)\bigr) = \varepsilon_i+\varepsilon_j+\varepsilon_k$. 
It follows that for vector fields $X = \sum_i X^i\cdot \partial_{x^i}$ and $Y = \sum_i Y^i \cdot \partial_{x^i}$, we have
\begin{moneq}
\nabla_XY = \sum_{ij} X^j\cdot \fracp{Y^i}{x^j} \cdot \partial_{x^i} 
+
\sum_{ijk} X^j \cdot \conjug^{\varepsilon_j}(Y^k)\cdot \Gamma_{jk}^i \cdot \partial_{x^i}
\end{moneq}%
%where $\mathfrak{C}$ denotes the conjugation in the odd variable on $\CA$ (by definition, $\mathfrak{C}(a_0+a_1) = a_0-a_1$ if $a_i \in \CA_i$). 
When the vector field $X$ is even, we have $\varepsilon(X^j) = \varepsilon_j$ and in that case the above formula can be written without signs as}
\begin{moneq}
\nabla_XY = \sum_{ij} X^j\cdot \fracp{Y^i}{x^j} \cdot \partial_{x^i}
+
\sum_{ijk} Y^k \cdot X^j \cdot \Gamma_{jk}^i \cdot \partial_{x^i}
\end{moneq}%
}

\begin{proclaim}{Corollary}
If $\nabla$ and $\nablah$ are connections on $M$ with Christoffel symbols $\Gamma_{jk}^i$ and $\Gammah_{jk}^i$ respectively, the tensor $S$ reads locally as
\begin{moneq}
S = 
\sum_{ijk}\extder x^k \otimes \extder x^j \cdot \Bigl( \Gamma_{jk}^i - \Gammah_{jk}^i \Bigr) \otimes \partial_{x^i}
\end{moneq}%
while the tensor $T$ is given by
\begin{align*}
T 
&= 
\sum_{ijk} \extder x^k \wedge \extder x^j \cdot \Gamma_{jk}^i(x)  \otimes \partial_{x^i}
\\&
=
\tfrac12\cdot \sum_{ijk} \extder x^k \wedge \extder x^j \cdot \bigl( \,\Gamma_{jk}^i - (-1)^{\varepsilon_j\varepsilon_k} \cdot \Gamma_{kj}^i \,\bigr)  \otimes \partial_{x^i}
\end{align*}
In particular $\nabla$ is torsion-free if and only if the Christoffel symbols are graded symmetric in the lower indices, \ie,  $\Gamma_{jk}^i = (-1)^{\varepsilon_j\varepsilon_k} \cdot \Gamma_{kj}^i$.

\end{proclaim}

If $y=(y^1, \dots, y^n)$ is another local system of coordinates, we can consider the Christoffel symbols $\Gammat_{jk}^i$ in terms of these coordinates:
\begin{moneq}
\Gammat_{jk}^i(y) 
= \contrf{\nabla_{\partial_{y^j}}\partial_{y^k}}{\,\extder y^i\vert_y}
\end{moneq}%
Now let $m\in M$ be the point in $M$ whose coordinates are $x$ or $y$ depending upon the choice of local coordinate system.
As tangent vectors transform as $\partial_{x^i}\vert_m = \sum_p (\partial_{x^i}y^p)(x)\cdot \partial_{y^p}\vert_m$, {it follows that the relation between $\Gamma$ and $\Gammat$ is given by}
%ON PEUT SUPPRIMER CE CALCUL ET DONNER DIRECTEMENT LA FORMULE \recalle{transChristoffel} ?
%
%\begin{align*}
%\sum_i \Gamma_{jk}^i(x) \cdot \partial_{x^i}\vert_m
%&=
%\Bigl(\ \nabla_{\partial_{x^j}}\bigl(\ \sum_r (\partial_{x^k}y^r)(x)\cdot \partial_{y^r}\, \bigr)\,\Bigr) \vert_m
%\\
%&
%=
%\sum_r (\partial_{x^j}\partial_{x^k}y^r)(x)\cdot \partial_{y^r}\vert_m
%\\&\qquad
%+
%\sum_r
%(-1)^{\varepsilon_j(\varepsilon_r+\varepsilon_k)} \cdot (\partial_{x^k}y^r)(x) \cdot \bigl(\nabla_{\partial_{x^j}} \partial_{y^r}\bigr)\vert_m
%\\
%&
%=
%\sum_p (\partial_{x^j}\partial_{x^k}y^p)(x)\cdot \partial_{y^p}\vert_m
%\\&\qquad
%+
%\sum_{qr}
%(-1)^{\varepsilon_j(\varepsilon_r+\varepsilon_k)} \cdot (\partial_{x^k}y^r)(x) \cdot (\partial_{x^j}y^q)(x) \cdot \bigl(\nabla_{\partial_{y^q}} \partial_{y^r} \bigr)\vert_m
%\\
%&
%=
%\sum_p (\partial_{x^j}\partial_{x^k}y^p)(x)\cdot \partial_{y^p}\vert_m
%\\&\qquad
%+
%\sum_{qr}
%(-1)^{\varepsilon_j(\varepsilon_r+\varepsilon_k)} \cdot (\partial_{x^k}y^r)(x) \cdot (\partial_{x^j}y^q)(x) \cdot \Gammat_q{}^p{}_r(y) \cdot \partial_{y^p}\vert_m
%\end{align*}
%which gives us the relations
\begin{multline}
\label{transChristoffel}
\quad
\sum_i \Gamma_{jk}^i(x) \cdot (\partial_{x^i}y^r)(x)
\\
=
(\partial_{x^j}\partial_{x^k}y^r)(x)
+
\sum_{s,t}
(-1)^{\varepsilon_j(\varepsilon_t+\varepsilon_k)} \cdot (\partial_{x^k}y^t)(x) \cdot (\partial_{x^j}y^s)(x) \cdot \Gammat_{st}^r(y)
\quad
\end{multline}

{Finally, let us consider $TM^{(0)}$ (the even part of the tangent bundle). With any local system of coordinates $x=(x^1, \dots, x^n)$ (resp. $y=(y^1, \dots, y^n)$) we associate the natural local system of coordinates $(x,v)$ (resp. $(y,w)$) on $TM^{(0)}$. More precisely, if $x$ are the coordinates of a point $m\in M$, then $(x,v)$ are the coordinates of the tangent vector $\tanvec = \sum_i v^i \cdot \partial_{x^i}\vert_m \in T_mM^{(0)}$. Now if $(x,v)$ and $(y,w)$ are the local coordinates of the same tangent vector $\tanvec$, \ie,
\begin{moneq}
\tanvec = \sum_i v^i \cdot \partial_{x^i}\vert_m = \sum_p w^p \cdot \partial_{y^p}\vert_m
\end{moneq}%
then we have}
\begin{moneq}[changeoftangcoord]
w^p = \sum_i v^i\cdot (\partial_{x^i}y^p)(x)
\end{moneq}%
{It follows that we have}
\begin{subequations}
\label{transeventangent}
\begin{align}
\label{transeventangentfirst}
\partial_{x^i}\vert_{\tanvec}
&=
\sum_p (\partial_{x^i}y^p)(x) \cdot \partial_{y^p}\vert_{\tanvec}
+
\sum_{jp} (-1)^{\varepsilon_i\varepsilon_j} v^j \cdot (\partial_{x^i}\partial_{x^j} y^p)(x) \cdot \partial_{w^p} \vert_{\tanvec}
\\
\label{transeventangentsecond}
\partial_{v^i}\vert_{\tanvec} 
&= 
\sum_p (\partial_{x^i}y^p)(x) \cdot \partial_{w^p}\vert_{\tanvec}
\end{align}
\end{subequations}

{With these preparations at hand, we now attack the question of defining geodesics.}
 We start very na\"ively in local coordinates and copy the classical case: a geodesic is a map $\gamma : \CA_0\to M$ given in local coordinates by $\gamma(t) = (\gamma^1(t), \dots, \gamma^n(t))$ satisfying the equations
\begin{moneq}[naivegeodesicsequations]
\frac{\extder^2\gamma^i}{\extder t^2}(t)
=
-
\sum_{jk}
\frac{\extder\gamma^k}{\extder t}(t) \cdot \frac{\extder\gamma^j}{\extder t}(t) \cdot
\Gamma_{jk}^i(\gamma(t)) 
\end{moneq}%
But to solve second order differential equations one needs initial conditions, which in our case are a starting point $x$ and an initial velocity $v$. And then the geodesic $\gamma$ depends upon these initial conditions, forcing us to write $\gamma_{(x,v)}$ instead of simply $\gamma$ and adding the initial conditions 
\begin{moneq}
\gamma^i_{(x,v)}(0) = x^i
\qquad\text{and}\qquad
\frac{\extder\gamma^i_{(x,v)}}{\extder t}(0)
=
v^i
\end{moneq}
It is here that our definition deviates from the one given in \cite{GW}, as we look at maps defined on $\CA_0\times TM^{(0)}$ rather than on $\CA_0\times \CA_1$ or an arbitrary product $\CA_0\times S$.
We now recall that any system of second order differential equations on a manifold can be expressed as a system of first order differential equations on the tangent bundle. This means that we look at curves $\gammat_{(x,v)}:\CA_0\to TM^{(0)}$ given in local coordinates by
\begin{moneq}
\gammat_{(x,v)}(t) = (\gamma^1_{(x,v)}(t), \dots, \gamma^n_{(x,v)}(t), \gammab^1_{(x,v)}(t), \dots, \gammab^n_{(x,v)}(t))
\end{moneq}%
satisfying the equations
\[
\left\{
\begin{array}{rcl}
\frac{\extder \gamma^i_{(x,v)}}{\extder t}(t)
&=&
\gammab^i_{(x,v)}(t)
\\[2\jot]
\frac{\extder \gammab^i_{(x,v)}}{\extder t}(t)
&=&
- \sum_{jk}\gammab^k_{(x,v)}(t) \cdot \gammab^j_{(x,v)}(t) \cdot \Gamma_{jk}^i(\gamma(t)) 
\end{array}
\right.
\]
and with initial conditions
\begin{moneq}
\gamma^i_{(x,v)}(0)
=
x^i
\qquad\text{and}\qquad
\gammab^i_{(x,v)}(0)
=
v^i
\end{moneq}%
We now recognize that these are exactly the equations of the integral curves of a vector field on $TM^{(0)}$. And indeed, using the Christoffel symbols we can define a vector field $G$ on $TM^{(0)}$ in local coordinates $(x,v)$ by
\begin{moneq}[definitionofG]
G\vert_{\tanvec}
=
\sum_i v^i \partial_{x^i}\vert_{\tanvec}
-
\sum_{ijk}
v^k\cdot v^j\cdot\Gamma_{jk}^i(x) 
\cdot \partial_{v^i}\vert_{\tanvec}
\end{moneq}%
Combining \recalle{transChristoffel} and \recalle{transeventangent}, it is immediate that these local expressions glue together to form a well defined global vector field $G$ on $TM^{(0)}$. As it is an even vector field, it has a flow $\Psi$ defined in an open subset $W_G$ of $\CA_0 \times TM^{(0)}$ containing $\{0\}\times TM^{(0)}$ and with values in $TM^{(0)}$ \cite[V.4.9]{Tu1}. In local coordinates we will write $\Psi(t,x,v) = (\Psi_1(t,x,v), \Psi_2(t,x,v))$, where $\Psi_1=(\Psi_1^1, \dots, \Psi_1^n)$ represents the base point while $\Psi_2=(\Psi_2^1, \dots, \Psi_2^n)$ represents the tangent vector. By definition of a flow, these functions thus satisfy the equations
\[
\left\{
\begin{array}{rcl}
\fracp{\Psi_1^i}{t}(t,x,v)
&=&
\Psi_2^i(t,x,v)
\\[2\jot]
\fracp{\Psi_2^i}{t}(t,x,v)
&=&
-
\sum_{jk}
\Psi_2^k(t,x,v) \cdot \Psi_2^j(t,x,v) \cdot
\Gamma_{jk}^i(\Psi_1(t,x,v)) 
\end{array}
\right.
\]
together with the initial conditions
\begin{moneq}
\Psi_1(0,x,v) = x
\qquad\text{and}\qquad
\Psi_2(0,x,v) = v
\end{moneq}%
{With the global vector field $G$ we thus have found an intrinsic coordinate free description of the equations we wrote for the geodesic curves $\gammat_{(x,v)}(t)$ and we are now in position to state a definition.}

\begin{definition}{Definition}
Let $\nabla$ be a connection in $TM$, let $\pi:TM^{(0)}\to M$ denote the canonical projection, let $G$ be the even vector field \recalle{definitionofG} and let $\Psi:W_G\to TM^{(0)}$ be its flow. For a fixed $(x,v)\cong\tanvec\in TM^{(0)}$ we will call the map $\gamma:\CA_0\to M$ defined by
\begin{moneq}
\gamma(t) = \pi\bigl(\Psi(t,\tanvec)\bigr) \cong \Psi_1(t,x,v)
\end{moneq}%
{\em the geodesic through $x\in M$ with initial velocity $v$.} Note that if $\tanvec$ is not in the body of $TM^{(0)}$, this curve is not (necessarily) smooth (see \cite[III.1.23g, V.3.19]{Tu1}).
\end{definition}

\begin{definition}{Remark}
One could define a similar vector field on $TM^{(1)}$, the odd part of the tangent bundle. More precisely, we denote by $(x,\vb)$ local coordinates on $TM^{(1)}$, where $(x,\vb)$ represents the tangent vector $\mathcal{V} = \sum_i \vb^i\cdot \partial_{x^i}\vert_m$, but the parity of $\vb^i$ is reversed: $\varepsilon(\vb^i) = \varepsilon_i+1 = \varepsilon(x^i)+1$. It thus is an odd tangent vector. These coordinates still change according to \recalle{changeoftangcoord} (with $v$ replaced by $\vb$), but an additional sign appears in the transformation of the tangent vectors: \recalle{transeventangentfirst} is replaced by
\begin{subequations}
\label{transeventangentodd}
\begin{align}
\label{transeventangentoddfirst}
\partial_{x^i}\vert_{\tanvec}
&=
\sum_p (\partial_{x^i}y^p)(x) \cdot \partial_{y^p}\vert_{\tanvec}
+
\sum_{jp} (-1)^{\varepsilon_i(\varepsilon_j+1)} \vb^j \cdot (\partial_{x^i}\partial_{x^j} y^p)(x) \cdot \partial_{\overline w^p} \vert_{\tanvec}
\end{align}
\end{subequations}
The analogon of the vector field $G$ on $TM^{(0)}$ would be the odd vector field $G'$ on $TM^{(1)}$ defined in local coordinates as
\begin{moneq}
G'\vert_{\tanvec}
=
\sum_i \vb^i \partial_{x^i}\vert_{\tanvec}
-
\sum_{ijk}
(-1)^{\varepsilon_k} \cdot \vb^k\cdot \vb^j\cdot\Gamma_{jk}^i(x) 
\cdot \partial_{\vb^i}\vert_{\tanvec}
\end{moneq}%
The transformation properties \recalle{transChristoffel}, \recalle{transeventangentsecond} and \recalle{transeventangentoddfirst} ensure that $G'$ is a well defined global vector field. However, the condition for an odd vector field to be integrable (with an odd time parameter $\tau$) is that its auto-commutator is zero \cite[V.4.17]{Tu1}. But the auto-commutator $[G',G']$ is given by 
\begin{align*}
{[G',G']} 
&=
-2\cdot \sum_{ijk} (-1)^{\varepsilon_k}\cdot \vb^k \cdot \vb^j \cdot \Gamma_{jk}^i(x) \cdot \partial_{x^i}
+
\text{ terms in $\partial_{\vb^i}$}
\\&
=
-\sum_{ijk} (-1)^{\varepsilon_k}\cdot \vb^k \cdot \vb^j \cdot (\Gamma_{jk}^i(x) - (-1)^{\varepsilon_j\varepsilon_k}\cdot \Gamma_{jk}^i(x)) \cdot \partial_{x^i}
+
\text{ terms in $\partial_{\vb^i}$}
\end{align*}
If this is to be zero, then at least the coefficients of $\partial_{x^i}$ have to be zero. But this is the case if and only if the connection $\nabla$ is torsion-free (on the odd tangent bundle, the combination $(-1)^{\varepsilon_k} \cdot \vb^k \cdot \vb^j$ is graded anti-symmetric). Moreover, if this is the case, then the vector field $G'$ reduces to $G' = \sum_i \vb^i \partial_{x^i}$, of which the auto-commutator indeed is zero (hence we don't have to compute the coefficients of $\partial_{\vb^i}$). But for this vector field the flow $\Phi'$ is given by:
\begin{moneq}
\Phi'(\tau,x,\vb) = (x+\tau\cdot \vb, \vb)
\end{moneq}%
which is rather uninteresting: the ``odd geodesics'' are ``straight odd lines'' in the direction of the tangent vector. Another way to see that this must happen is {the following set of observations. If we use an odd time parameter $\tau$, it follows immediately that the velocity vector should be an odd tangent vector. Moreover, when we write the naïve equations \recalle{naivegeodesicsequations} for the geodesics, the left hand side is identically zero because $\partial_\tau \scirc \partial_\tau = 0$. And then this equation tells us that the connection should be torsion-free. We are thus left with the condition that the connection should be torsion-free, together with the initial conditions $\gamma(0,x,\vb) = x$ and $\partial_\tau\gamma(0,x,\vb) = \vb$. And these give us our straight odd lines.}

\end{definition}

\section{Projective equivalence}\label{proj}

We now consider the situation in which we have two connections $\nabla, \nablah$ on $M$ and we wonder under what conditions these two connections have ``the same'' geodesics as images in $M$. More precisely, if $\Psi(t,\tanvec)$ and $\Psih(t,\tanvec)$ are the geodesic flows for $\nabla$ and $\hat{\nabla}$ respectively, the na\"ive question is under what conditions we have 
\begin{moneq}
\{ \, \Psi_1(t,\tanvec) :  t\in \CA_0 \,\}
=
\{ \, \Psih_1(t,\tanvec) :  t\in \CA_0 \,\}
\end{moneq}% 
A more precise question is under what conditions we can find a reparametrization function $r:\CA_0 \times TM \to \CA_0$ such that we have
\begin{moneq}[reparamPsi1]
\forall t\in \CA_0\quad:\quad
\Psi_1(r(t,\tanvec), \tanvec) = \Psih_1(t,\tanvec)
\end{moneq}%
Note that we added an explicit dependence on the initial condition $\tanvec$ in the repara\-me\-trization function $r$, as there is no reason that geodesics through different points should be reparametrized in the same way.

\begin{definition}{Definition}
We say that the connections $\nabla$ and $\nablah$ have \emph{the same  geodesics up to reparametrization} if there exists a function $r : \CA_0 \times TM \to \CA_0$ such that $r(0,\tanvec) = 0$, $(\partial r/\partial t)(0,\tanvec) = 1$ and for which equation (\ref{reparamPsi1}) holds.\footnote{The additional conditions $r(0,\tanvec) = 0$ and $(\partial r/\partial t)(0,\tanvec) = 1$ ensure that the reparametrization transforms each geodesic of $\nabla$ into the geodesic of $\nablah$ with the same initial conditions.}
\end{definition}

We are going to characterize the connections that %\textcolor{blue}{j'ai remplacé ``which'' par ``that'', car mes souvenirs de grammaire anglaise me disent que ``which'' est descriptif et ``that'' est restrictif} %
have the same geodesics up to reparametrization in terms of the form of the tensor $S$ which measures the difference between these two connections. In order to do that, we are going to proceed in two steps. First, we show that (\ref{reparamPsi1}) holds if and only if the geodesic flow $\Psi$ of $G$, the (difference) tensor $S=\nabla - \nablah$ and the reparametrization function $r$ are related through a certain differential equation.

\begin{proclaim}{Proposition}
The connections $\nabla$ and $\hat{\nabla}$ have the same geodesics up to reparametrization if and only if there exists a function $r:\CA_0 \times TM \to \CA_0$ such that $r(0,\tanvec) = 0$, $(\partial r/\partial t)(0,\tanvec) = 1$ and for which the following differential equation holds:
\begin{multline}
\label{reparametrization}
\frac{\partial^2r}{\partial t^2}(t,\tanvec) \cdot\frac{\partial\Psi_1}{\partial t}(r(t,\tanvec),\tanvec) 
\\
= 
\Bigl(\fracp{r}{t}(t,\tanvec)\Bigr)^2 \cdot
S_{\Psi_1(r(t,\tanvec),\tanvec)}\Bigl(\,\frac{\partial\Psi_1}{\partial t}(r(t,\tanvec),\tanvec)\,,\, \frac{\partial\Psi_1}{\partial t}(r(t,\tanvec),\tanvec) \,\Bigr)
\end{multline}
\end{proclaim}
\begin{preuve}
Let us show that the condition is necessary. In view of (\ref{naivegeodesicsequations}), if $\Psi_1(r(t,\tanvec), \tanvec)$ is a geodesic for $\hat{\nabla}$, then 
\begin{moneq}
0 = \frac{\partial^2\Psi_1^i(r(t,\tanvec),\tanvec)}{\partial t^2} 
+ \sum_{j,k}
\frac{\partial\Psi_1^k(r(t,\tanvec),\tanvec)}{\partial t} \cdot \frac{\partial\Psi_1^j(r(t,\tanvec),\tanvec)}{\partial t}\cdot
\hat{\Gamma}_{jk}^i(\Psi_1(r(t,\tanvec), \tanvec))
\end{moneq}%
Let us replace in this equation $\hat{\Gamma}_{jk}^{i}$ by $\Gamma_{jk}^{i}-S_{jk}^{i}$ and let us apply the chain rule to compute the derivatives of the functions $\Psi_1^i(r(t,\tanvec),\tanvec)$. Doing so, we obtain %\textcolor{blue}{\bf j'ai enlevé du texte}%
%that if $\Psi_1(r(t,\tanvec), \tanvec)$ is a geodesic for $\hat{\nabla}$, then
\begin{eqnarray*}
0 
&=& 
\frac{\partial^2r}{\partial t^2}(t,\tanvec) \cdot\frac{\partial\Psi_1}{\partial t}(r(t,\tanvec),\tanvec)
+ \left(\frac {\partial r}{\partial t}(t,\tanvec)\right)^{2} \left(\frac{\partial^2\Psi_1^i}{\partial t^2}(r(t,\tanvec),\tanvec)\right) 
\\
&&+ \left(\frac {\partial r}{\partial t}(t,\tanvec)\right)^{2} \left(\sum_{j,k}
\frac{\partial\Psi_1^k}{\partial t}(r(t,\tanvec),\tanvec) \cdot \frac{\partial\Psi_1^j}{\partial t}(r(t,\tanvec),\tanvec)\cdot\Gamma_{jk}^i(\Psi_1(r(t,\tanvec), \tanvec))\right)
\\
&&-\Bigl(\fracp{r}{t}(t,\tanvec)\Bigr)^2 \left( \sum_{j,k}
\frac{\partial\Psi_1^k}{\partial t}(r(t,\tanvec),\tanvec) \cdot \frac{\partial\Psi_1^j}{\partial t}(r(t,\tanvec),\tanvec)\cdot S_{jk}^i(\Psi_1(r(t,\tanvec), \tanvec))\right) 
\end{eqnarray*}
Using the fact that $\Psi_1$ is a geodesic for $\nabla$, 
%\[\frac{\partial^2\Psi_1^i}{\partial t^2}(r(t,\tanvec),\tanvec)+\sum_{jk}
%\frac{\partial\Psi_1^k}{\partial t}(r(t,\tanvec),\tanvec) \cdot \frac{\partial\Psi_1^j}{\partial t}(r(t,\tanvec),\tanvec)\cdot\Gamma_{jk}^i(\Psi_1(r(t,\tanvec), \tanvec))=0,\]
the second and third term on the right hand side cancel and hence %we eventually see that 
this equation reduces to (\ref{reparametrization}).

In order to show the converse, it suffices to note that the above computations also show that if \eqref{reparametrization} is satisfied, then the curve 
%thanks to the computations above that the curve 
\[\left(\Psi_1(r(t,\tanvec), \tanvec),\frac {\partial r}{\partial t}(t,\tanvec)\cdot\frac{\partial\Psi_1}{\partial t}(r(t,\tanvec),\tanvec)\right)
\]
satisfies the equation of the flow $(\Psih_{1}(t,\tanvec),\Psih_{2}(t,\tanvec))$ of $\hat{G}$, the geodesic vector field corresponding to $\hat{\nabla}$. As it satisfies the same initial conditions as $(\Psih_{1}(t,\tanvec),\Psih_{2}(t,\tanvec))$ at $t=0$, these two curves have to coincide, and in particular %. In particular, we obtain that 
$\Psi_1(r(t,\tanvec), \tanvec) = \Psih_1(t,\tanvec)$.   
\end{preuve}

Now in order to obtain Weyl's characterization in the super context, it remains to show that condition (\ref{reparametrization}) amounts to imposing that $S$ can be expressed by means of an even (super) $1$-form. As for the previous Proposition, the proof of the theorem follows the lines of the classical case. It invokes a technical Lemma which roughly says that if we have a bilinear function $S(v,w)$ such that $S(v,v)=h(v)\cdot v$ for some function $h$, then $h$ must be linear in $v$. The proof of this technical Lemma is elementary but long, simply because we have to be careful with the odd coordinates and moreover, everything depends upon additional parameters (the local coordinates $x$ and $\xi$ on $M$). Therefore the proof of the lemma will be given after that of the Theorem.
%\textcolor{blue}{J'ai changé l'ordre, car c'est dans le lemme qu'il faut faire attention. Et j'ai changé quadratique en bilinéaire, sinon on ne peut pas écrire $S(v,v)$ comme fonction de DEUX variables.}% 

%, but here we have to be careful because we have odd coordinates and moreover, everything depends upon additional parameters (the $x$ and $\xi$). In practice, the proof of the theorem will invoke a technical Lemma to show that if $S$ is quadratic in $v$ and satisifies $S(v,v) = h(v) \cdot v$ for a certain function $h$, then this $h$ must be linear in $v$. The proof of the Lemma is elementary but lengthy. Therefore, this proof will be given after that of the theorem.

\begin{proclaim}[linearizationlemma]{Lemma}
Let $E$ be a graded vector space of graded dimension $p\vert q$ with even basis vectors $e_1, \dots, e_p$ and odd basis vectors $f_1, \dots, f_q$, let $U$ be an open coordinate subset of a manifold $M$ with local even coordinates $x$ and local odd coordinates $\xi$. Suppose that $S:U\times E\times E \to E$ is a smooth function which is left-bilinear, graded symmetric in the product $E\times E$ and for which there is a smooth function $h:U\times E_0\to \CA$ such that
\begin{moneq}[determiningoneform]
\forall (x,\xi)\in U\ \forall v\in E_0 : 
S(x,\xi,v,v) = h(x,\xi,v) \cdot v
\end{moneq}%
Then there exists a unique smooth function $\alpha:U\to E^*$ such that $h(x,\xi,v) = \contrf{v}{\alpha(x,\xi)}$ and 
\begin{moneq}
S(x,\xi,v,w) = \tfrac12 \cdot \bigl(\, \contrf{v}{\alpha(x,\xi)}\cdot w + \contrf{w}{\alpha(x,\xi)}\cdot v \,\bigr)
\end{moneq}%
\end{proclaim}

\begin{proclaim}[conformalequivalenceconnections]{Theorem}
Two torsion-free connections $\nabla$ and $\nablah$ on $M$ have the same geodesics up to reparametrization if and only if there exists a smooth even $1$-form $\alpha$ on $M$ such that the tensor $S = \nabla - \nablah$ is given by 
\begin{moneq}[Sintermsofalpha]
S_x(v,w) = \tfrac12\cdot ( \contrf{v}\alpha_x \cdot w + (-1)^{\varepsilon(v) \cdot \varepsilon(w)} \cdot \contrf w\alpha_x \cdot v )
\end{moneq}%
for any $x \in M$ and any homogeneous $v,w \in T_xM$.
\end{proclaim}

\begin{preuve}[Proof of the theorem]%[Proof of \recallp{conformalequivalenceconnections}]
We first assume that we have a reparametrization $r$ that transforms the geodesics of $\nabla$ into those of $\nablah$. Taking $t=0$ in \recalle{reparametrization} and using the initial conditions for $\Psi$ and $r$, we get the following (vector) equation in local coordinates:
\begin{moneq}[reduceddefeqforreparametrization]
v \cdot \frac{\partial^2 r}{\partial t^2}(0,x,v)
=
S_x(v,v)
\end{moneq}%
Lemma \ref{linearizationlemma}, with $h$ being here the function $h(x,v) = \frac{\partial^2r}{\partial t^2}(0,x,v)$, gives us a (local) smooth $1$-form $\alpha$, which must be even by parity considerations. But \recalle{reduceddefeqforreparametrization} is an intrinsic equation which does not depend upon the choice of local coordinates (because \recalle{reparametrization} is intrinsic). As the $1$-form $\alpha$ is unique, %Also the $1$-form $\alpha$ is unique, hence
the local $1$-forms $\alpha$ given by Lemma \ref{linearizationlemma} glue together to form a global smooth even $1$-form $\alpha$ satisfying \recalle{Sintermsofalpha}.

To show the converse, let us now assume that we have an even $1$-form $\alpha$ on $M$ such that the tensor $S$ is given by \recalle{Sintermsofalpha}. Then \recalle{reparametrization} reduces to the (vector) equation
\begin{multline}
\label{defeqforreparametrizationgivenS}
\frac{\partial^2 r}{\partial t^2}(t,x,v) \cdot \frac{\partial\Psi_1}{\partial t}(r(t,x,v),x,v) 
\\
=
\Bigl( \fracp rt(t,x,v)\Bigr)^2 \cdot
\iota\left( \frac{\partial\Psi_1}{\partial t}(r(t,x,v),x,v)\right){\alpha_{\Psi_1(r(t,x,v),x,v)}} \cdot \frac{\partial\Psi_1}{\partial t}(r(t,x,v),x,v) 
\end{multline}
For this to be true for all geodesics of $\nabla$, the function $r$ thus has to satisfy the second order differential equation
\begin{align*}
\frac{\partial^2 r}{\partial t^2}(t,x,v)
&=
\Bigl( \fracp rt(t,x,v)\Bigr)^2 \cdot \iota \left(\frac{\partial\Psi_1}{\partial t}(r(t,x,v),x,v)\right){\alpha_{\Psi_1(r(t,x,v),x,v)}}
\end{align*}
As for the geodesic equations, we translate this into a system of first order differential equations by introducing a second function $s:\CA_0\times TM^{(0)}\to \CA_0$ and we obtain
\[
\left\{
\begin{array}{rcl}
\fracp{r}{t}(t,x,v)
&=&
s(t,x,v)
\\[2\jot]
\fracp{s}{t}(t,x,v)
&=&
s(t,x,v)^2 \cdot \iota \left(\frac{\partial\Psi_1}{\partial t}(r(t,x,v),x,v)\right){\alpha_{\Psi_1(r(t,x,v),x,v)}}
\end{array}
\right.
\]
while the initial conditions for $r$ yield $r(0,x,v) = 0$ and $s(0,x,v) = 1$. To show that these equations always have a (unique) solution, we just note that these equations determine the flow of the even vector field $R$ on $(\CA_0)^2 \times TM^{(0)}$ given by
\begin{moneq}
R\vert_{(r,s,\tanvec)}
=
s\cdot \fracp{}{r}
+
s^2 \cdot \iota\left({\frac{\partial\Psi_1}{\partial t}(r,\tanvec)}\right){\alpha_{\Psi_1(r,\tanvec)}} \cdot \fracp{}{s}
\end{moneq}%
And indeed, the equations for the flow $\Phi= (\Phi_r, \Phi_s, \Phi_1, \Phi_2)$ of $R$ are given by
\[
\left\{
\begin{array}{rcl}
\fracp{\Phi_r}{t}(t,r_o,s_o,x,v)
&=&
\Phi_s(t,r_o,s_o,x,v)
\\[2\jot]
\fracp{\Phi_s}{t}(t,r_o,s_o,x,v)
&=&
(\Phi_s(t,r_o,s_o,x,v))^2 
\\&&\qquad
\cdot \iota\left({\frac{\partial\Psi_1}{\partial t}(\Phi_r(t,r_o,s_o,x,v),x,v)}\right) {\alpha_{\Psi_1(\Phi_r(t,r_o,s_o,x,v),x,v)}} 
\\[2\jot]
\fracp{\Phi_1}{t}(t,r_o,s_o,x,v)
&=&
0
\\[2\jot]
\fracp{\Phi_2}{t}(t,r_o,s_o,x,v)
&=&
0
\end{array}
\right.
\] 
It thus suffices to define $r(t,\tanvec) = \Phi_r(t,0,1,\tanvec)$ and $s(t,\tanvec) = \Phi_s(t,0,1,\tanvec)$ to obtain the desired functions. 
\end{preuve}

%Before we attack the proof of this theorem, we fist prove a preliminary lemma. Its proof is an elementary computation, which is lengthy only because we have to be careful with odd coordinates. The content of this lemma is a rather easy statement: if we have a function $S$ which is quadratic in $v$ and which is equal to $v$ times a function $h$ of $v$, then $h$ must be linear in $v$. This issue is obscured by the fact that we have a vector valued function, that we have odd coordinates and that everything depends upon additional parameters (the $x$ and $\xi$).

\begin{preuve}[Proof of the lemma]
Uniqueness of $\alpha$ follows from the equation $h(x,\xi,v) = \contrf{v}{\alpha(x,\xi)}$. To prove existence, let us start by introducing global (linear, left) coordinates $y,\eta$ on $E_0$ by
\begin{moneq}
v\in E_0 \quad\Rightarrow\quad
v=\sum_i y_i \cdot e_i + \sum_i \eta_i \cdot f_i
\end{moneq}%
Using bilinearity and graded symmetry, we thus can write
\begin{multline*}
S(x,\xi,v,v)
=
\sum_{i,j} y_i y_j \cdot S(x,\xi,e_i,e_j)
\\
+
2\sum_{i,j} y_i \eta_j \cdot S(x,\xi,e_i,f_j)
+
\sum_{i,j} \eta_j\eta_i \cdot S(x,\xi,f_i,f_j)
\end{multline*}
The functions $S$, when evaluated in a pair of basis vectors of $E$, is a smooth function on $U$ with values in $E$. As such we can determine the coefficients with respect to the given basis for $E$ as for instance
\begin{moneq}
S(x,\xi,e_i,e_j)
=
\sum_p S_p(x,\xi,e_i,e_j) \cdot e_p 
+
\sum_p \sigma_p(x,\xi,e_i,e_j) \cdot f_p
\end{moneq}%
When we substitute this in \recalle{determiningoneform} with the (linear, left) coordinates of $v\in E_0$, we get the system of equations
\begin{align*}
h(x,\xi,y,\eta) \cdot y_p
&=
\sum_{i,j} y_i y_j \cdot S_p(x,\xi,e_i,e_j)
\\&
\qquad+
2\sum_{i,j} y_i \eta_j \cdot S_p(x,\xi,e_i,f_j)
+
\sum_{i,j} \eta_j\eta_i \cdot S_p(x,\xi,f_i,f_j)
\\
h(x,\xi,y,\eta) \cdot \eta_p
&=
\sum_{i,j} y_i y_j \cdot \sigma_p(x,\xi,e_i,e_j)
\\&
\qquad+
2\sum_{i,j} y_i \eta_j \cdot \sigma_p(x,\xi,e_i,f_j)
+
\sum_{i,j} \eta_j\eta_i \cdot \sigma_p(x,\xi,f_i,f_j)
\end{align*}
%Note that we put the (linear) factor $\eta_p$ to the right of the function $h$, as the parity of $h$ is not given.
Applying \recallp{expandinginoddpowers} we can expand these equations in powers of the $\xi$ coordinates and equate the separate powers $\xi^J$ giving
\begin{align}
h_J(x,y,\eta) \cdot y_p
&=
\sum_{i,j} y_i y_j \cdot S_{p,J}(x,e_i,e_j)
\label{evenbasisequation}
+
2\sum_{i,j} y_i \eta_j \cdot (-1)^{\vert J\vert} \cdot S_{p,J}(x,e_i,f_j)
\\&\qquad
+
\sum_{i,j} \eta_j\eta_i \cdot S_{p,J}(x,f_i,f_j)
\notag
\\
h_J(x,y,\eta) \cdot \eta_p
&=
\sum_{i,j} y_i y_j \cdot \sigma_{p,J}(x,e_i,e_j)
\label{oddbasisequation}
+
2\sum_{i,j} y_i \eta_j \cdot (-1)^{\vert J\vert} \cdot \sigma_{p,J}(x,e_i,f_j)
\\&\qquad
+
\sum_{i,j} \eta_j\eta_i \cdot \sigma_{p,J}(x,f_i,f_j)
\notag
\end{align}
Note that we had to add a factor $(-1)^{\vert J\vert}$ in the right hand side for the terms linear in $\eta$, because we factor the powers of $\xi$ to the left, and interchanging a power $\xi^J$ with a linear factor $\eta$ gives this sign.
We now expand the functions $h_J$ in powers of the odd coordinates $\eta$:
\begin{multline*}
h_J(x,y,\eta) = h_{J,\emptyset}(x,y) + \sum_q \eta_q \cdot h_{J,\{q\}}(x,y) \\
+ \sum_{q<r} \eta_q\eta_r \cdot h_{J,\{q,r\}}(x,y) + \sum_{I, \vert I\vert \ge 3} \eta^I \cdot h_{J,I}(x,y)
\end{multline*}%
When we now invoke \recallp{expandinginoddpowers} applied to \recalle{evenbasisequation}, we get the equations
\begin{align}
h_{J,\emptyset}(x,y) \cdot y_p
&= 
\sum_{i,j} y_iy_j \cdot S_{p,J}(x,e_i,e_j) 
\label{zerothordereven}
\\
h_{J,\{q\}}(x,y) \cdot y_p
&=
2\sum_i y_i \cdot (-1)^{\vert J\vert} \cdot S_{p,J}(x,e_i,f_q)
\label{firstordereven}
\\
h_{J,\{q,r\}}(x,y) \cdot y_p
&=
2S_{p,J}(x,f_r,f_q)
\label{secondordereven}
\\
h_{J,I}(x,y) \cdot y_p
&=
0 \qquad \vert I\vert \ge 3
\label{thirdordereven}
\end{align}
As these are equations between smooth functions of even coordinates only, we may consider them to be equations of smooth functions of real coordinates. And remember, the $y$ coordinates run over the whole of $\RR$ as they are coordinates on a (graded) vector space. These functions thus are in particuler smooth at $y=0$.

As the right hand sides of \recalle{secondordereven} and \recalle{thirdordereven} do not depend upon the $y$ coordinates and their left hand sides have at least degree one in $y$, it follows that the coefficients must be zero, and thus the right hand side of \recalle{secondordereven} too:
\begin{moneq}
h_{J,I}(x,y) = 0 \text{ for } \vert I\vert\ge3
\quad,\quad
h_{J,\{q,r\}}(x,y) = 0
\quad,\quad
S_{p,J}(x,f_r,f_q) = 0
\end{moneq}%
From \recalle{firstordereven} it follows easily that $h_{J,\{q\}}(x,y)$ is independent of the $y$ coordinates:
\begin{moneq}
h_{J,\{q\}}(x,y) = h_{J,\{q\}}(x)
\end{moneq}% 
and that we must have
\begin{moneq}
(-1)^{\vert J\vert} \cdot S_{p,J}(x,e_i,f_q)
=
\tfrac12\cdot
\delta_{ip}\cdot h_{J,\{q\}}(x)
\end{moneq}%
Using the bilinearity of $S$, one can show that \recalle{zerothordereven} implies that $h_{J,\emptyset}(x,y)$ must be linear in $y$:
\begin{moneq}
h_{J,\emptyset}(x,y)
=
\sum_q h_{J,\emptyset}^q(x)\cdot y_q
\end{moneq}%
and then that we must have
\begin{moneq}
S_{p,J}(x,e_i,e_j)
=
\tfrac12\cdot \Bigl( \delta_{ip} \cdot h_{J,\emptyset}^j(x) + \delta_{jp} \cdot h_{J,\emptyset}^i(x) \Bigr)
\end{moneq}%

We now apply exactly the same reasoning to \recalle{oddbasisequation}, equating the separate powers of $\eta$ and using what we already know about the functions $h_{J,I}(x,y)$. This gives us the equations
\begin{align}
0
&=
\sum_{i,j} y_iy_j \cdot \sigma_{p,J}(x,e_i,e_j)
\label{sigmabothes}
\\
\sum_q y_q \cdot h_{J,\emptyset}^q(x)
&=
2\sum_i y_i \cdot (-1)^{\vert J\vert} \cdot \sigma_{p,J}(x,e_i,f_p)
\label{sigmaeandf}
\\
\tfrac12 \cdot \Bigl( h_{J,\{j\}}(x)\cdot \delta_{ip} - h_{J,\{i\}}(x)\cdot \delta_{jp} \Bigr)
&=
\sigma_{p,J}(x,f_i, f_j)
\end{align}
As these are (again) equations between smooth functions of real variables, we may conclude from \recalle{sigmabothes} that we have $\sigma_{p,J}(x,e_i,e_j)=0$ and from \recalle{sigmaeandf} that we have $\tfrac12 \cdot h_{J,\emptyset}^i(x) \cdot \delta_{jp} = (-1)^{\vert J\vert} \cdot \sigma_{p,J}(x,e_i,f_j)$.

To summarize, we have found the following equalities
\begin{align*}
h_J(x,y,\eta)
&=
\sum_q y_q \cdot h_{J,\emptyset}^q(x) + \sum_q \eta_q \cdot h_{J,\{q\}}(x)
\\
S_{p,J}(x,e_i,e_j)
&=
\tfrac12\cdot \Bigl( \delta_{ip} \cdot h_{J,\emptyset}^j(x) + \delta_{jp} \cdot h_{J,\emptyset}^i(x) \Bigr)
\\
S_{p,J}(x,e_i,f_j) 
&=
\tfrac12\cdot(-1)^{\vert J\vert} \cdot 
\delta_{ip}\cdot h_{J,\{j\}}(x)
\\
S_{p,J}(x,f_i,f_j)
&=
0
\\
\sigma_{p,J}(x,e_i,e_j)
&=
0
\\
\sigma_{p,J}(x,e_i,f_j) 
&=
\tfrac12 (-1)^{\vert J\vert} \cdot  \cdot h_{J,\emptyset}^i(x) \cdot \delta_{jp}
\\
\sigma_{p,J}(x,f_i,f_j)
&=
\tfrac12 \cdot \Bigl( h_{J,\{j\}}(x)\cdot \delta_{ip} - h_{J,\{i\}}(x)\cdot \delta_{jp} \Bigr)
\end{align*}

We now define the smooth functions $H_\emptyset^q, H_{\{q\}}:U\to \CA$ by
\begin{moneq}
H_\emptyset^q(x,\xi) 
=
\sum_J \xi^J \cdot h_{J,\emptyset}^q(x)
\quad,\quad
H_{\{q\}}(x,\xi)
=
\sum_J \xi^J \cdot (-1)^{\vert J\vert} \cdot h_{J,\{q\}}(x)
\end{moneq}%
Using these functions, we now put the powers of $\xi$ back in to obtain
{%\allowdisplaybreaks[1]
\begin{align*}
h(x, \xi ,y,\eta)
&=
\sum_J \xi^J \cdot \Biggl( \sum_q y_q \cdot h_{J,\emptyset}^q(x) + \sum_q \eta_q \cdot h_{J,\{q\}}(x) \Biggr)
\\&
=
\sum_q y_q \cdot \sum_J \xi^J \cdot h_{J,\emptyset}^q(x) + \sum_q \eta_q \cdot \sum_J \xi^J \cdot (-1)^{\vert J\vert} \cdot h_{J,\{q\}}(x)
\\&
=
\sum_q y_q \cdot H_{\emptyset}^q(x,\xi) + \sum_q \eta_q \cdot H_{\{q\}}(x,\xi)
\\
S_{p}(x,\xi,e_i,e_j)
&=
\sum_J \xi^J \cdot S_{p,J}(x,e_i,e_j)
%\\&
=
\tfrac12\cdot \Bigl( \delta_{ip} \cdot H_{\emptyset}^j(x,\xi) + \delta_{jp} \cdot H_{\emptyset}^i(x,\xi) \Bigr)
\\
S_{p}(x,\xi,e_i,f_j) 
&=
\sum_J \xi^J \cdot S_{p,J}(x,e_i,f_j) 
=
\tfrac12\cdot \delta_{ip}\cdot H_{\{j\}}(x,\xi)
\\
S_{p}(x,\xi,f_i,f_j)
&=
0
=
\sigma_{p}(x,\xi,e_i,e_j)
\\
\rho \cdot \sigma_{p}(x,\xi,e_i,f_j) 
&=
\rho \cdot \sum_J \xi^J \cdot \sigma_{p,J}(x,e_i,f_j) 
=
\tfrac12   H_{\emptyset}^i(x,\xi) \cdot \delta_{jp} \cdot \rho
\\
\rho \cdot \sigma_{p}(x,\xi,f_i,f_j)
&=
\rho \cdot \sum_J \xi^J \cdot \sigma_{p,J}(x,f_i,f_j)
\\&
=
\tfrac12 \cdot \Bigl( H_{\{j\}}(x,\xi)\cdot \delta_{ip} - H_{\{i\}}(x,\xi)\cdot \delta_{jp} \Bigr) \cdot \rho
\end{align*}}% end allowdisplaybreaks
where $\rho$ is any odd variable. Finally, we can reconstruct the full function $S$: if $v$ reads as $\sum_i y_i e_i + \sum_i \eta_i f_i$ and $w$ reads as $\sum_j z_j e_j, + \sum_j \zeta_j f_j$, then direct substitution gives us
\begin{align*}
S(x,\xi,v,w)
&=
\tfrac12 \Bigl(
\sum_j z_j \cdot H_{\emptyset}^j(x,\xi) + \sum_j \zeta_j H_{\{j\}}(x,\xi) 
\Bigr) \cdot v
\\&\qquad
+
\tfrac12 \Bigl(
\sum_j y_j \cdot H_{\emptyset}^j(x,\xi) + \sum_j \eta_j H_{\{j\}}(x,\xi) 
\Bigr) \cdot w
\end{align*}

This suggests that we introduce the left-linear form $\alpha:U\to E^*$ by 
\begin{moneq}
\contrf{v}{\alpha(x,\xi)}
=
\contrf{\sum_i y_i e_i + \sum_i \eta_i f_i}{\alpha(x,\xi)}
=
\sum_i y_i \cdot H_{\emptyset}^i(x,\xi) + \sum_i \eta_i \cdot H_{\{i\}}(x,\xi)
\end{moneq}%
where $y_i,\eta_i$ are arbitrary (non-homogeneous) coefficients. It then follows immediately that we have
\begin{moneq}
S(x,\xi,v,w)
=
\tfrac12 \cdot \bigl(\, \contrf{w}{\alpha(x,\xi)}\cdot v + \contrf{v}{\alpha(x,\xi)} \cdot w \,\bigr)
\end{moneq}%
It also follows that we have
\begin{moneq}
h(x,\xi,y,\eta)
=
\contrf{v}{\alpha(x,\xi)}
\end{moneq}%
confirming the equation $S(x,\xi,v,v) = h(x,\xi,v) \cdot v$ for even vectors $v$.
\end{preuve}

\section{Super metrics and connections}\label{Smetricconnections}

As in non-super geometry, connections on the tangent bundle arise naturally when the supermanifold is equipped with a metric. Moreover, again as in non-super geometry, geodesics in this context can be interpreted as the trajectories on the supermanifold of a free particle whose kinetic energy is given by the metric. We now substantiate these claims.
More precisely, we shall first expose some basic theory of super metrics and their associated Levi-Civita (super) connections. Then we shall briefly describe the mechanics of a free particle whose kinetic energy is given by the metric and finally, following \cite{GW}, we shall relate the Hamiltonian vector field of this mechanical system to the geodesic vector field of the corresponding \emph{metric connection}.

\begin{definition}{Definition}
A (super) \emph{metric} $g$ on a supermanifold $M$ is an even graded symmetric non-degenerate smooth section of the bundle $T^*M\otimes T^*M \to M$.
A \emph{Riemannian supermanifold} is a pair $(M,g)$ with $M$ a supermanifold and $g$ a metric on $M$.
\end{definition}

A metric $g$ on $M$ amounts to a collection of maps $g_m:T_mM\times T_mM \to \CA$ (depending smoothly on $m\in M$) possessing 
%\textcolor{red}{for all $m\in M$} 
the following four properties:
\begin{itemize}
\item
%For each $m \in M$, the 

The map $(v,w)\mapsto 
\contrf{v,w}{g_m}$ 
%\textcolor{green}{AU CHOIX} $(v,w)\mapsto g_m(v,w)$
is (left-)bilinear in $v$ and $w\,$;%
\footnote{Since the map $g_m$ is supposed to be even, we could also have written $g_m(v,w)$ instead of $\contrf{v,w}{g_m}$. However, once we express $g_m$ in terms of the left-dual basis $dx^i$, there is a high risk of confusion on how to compute evaluations, as we have $(dx^j)(\partial_{x^i}) = (-1)^{\varepsilon_{x^i}} \delta^j_i$, and not (as one might be inclined to think) $(dx^j)(\partial_i) = \delta^j_i$, simply because we have (by definition of the left-dual basis): $\delta^j_i=\contrf{\partial_{x^i}}{dx^j} = (-1)^{\varepsilon_i\varepsilon_j}\, (dx^j)(\partial_{x^i})$.}

\item
for all homogeneous $v,w\in T_mM$~: $\varepsilon(\contrf{v,w}{g_m}) = \varepsilon(v) + \varepsilon(w)$;
%\textcolor{green}{AU CHOIX} $\varepsilon(g_m(v,w)) = \varepsilon(v) + \varepsilon(w)$

\item
for all homogeneous $v,w\in T_m$~: 
%$g_m(w,v) = (-1)^{\varepsilon(v)\varepsilon(w)} \,g_m(v,w)$ \textcolor{green}{AU CHOIX} 
$\contrf{w,v}{g_m} = (-1)^{\varepsilon(v)\varepsilon(w)} \, \contrf{v,w}{g_m}$.

\end{itemize}
Now for each $m \in M$, the map $g_m$ can be seen as transforming tangent vectors into cotangent vectors, \ie, we can define a map $g_m^\flat : T_mM \to T_m^*M$ by setting
%\[
%g_m^\flat(v) = \contrf{v}g_m = g_m(\cdot,v)
%\]
%\textcolor{green}{AU CHOIX}
\begin{moneq}
\contrf{v}{g_m^\flat} = \contrf{v}{g_m} = \contrf{\cdot,v}{g_m}
\quad\text{\ie,}\quad
\contrf{w}{\bigl(\contrf{v}{g_m^\flat}\bigr)} = \contrf{w}{\bigl(\contrf{v}{g_m}\bigr)} 
\equiv
\contrf{w,v}{g_m}
\end{moneq}
With this definition we can state the the fourth condition 
\begin{itemize}
\item $g_m^\flat : T_mM \to T_m^*M$ is a (left-)linear bijection.
\end{itemize}
The collection of all maps $g_m^\flat$ gives rise to an even bundle isomorphism $g^\flat : TM \to T^*M$, whose inverse is denoted by $g^\sharp : T^*M \to TM$. As usual, the use of the musical superscripts is inspired by the fact that $g^\flat$ lowers indices of tensors, wheras $g^\sharp$ raises them.
%\textcolor{blue}{Votre façon de l'écrire suggérait que la non-dégénérescence est une conséquence des trois premières conditions. Et je n'ose pas écrire la condition de non-dégénérescence comme pour tout $v\neq0$ il existe $w$ tel que $g_m(v,w)\neq0$ à cause des nilpotents.}%

%For each $m \in M$, $g_m^\flat$ is a (left-)linear bijection an that the collection of those maps gives rise to an even bundle isomorphism $g^\flat : TM \to T^*M$, whose inverse is denoted by $g^\sharp : T^*M \to TM$. As usual, the use of the musical superscripts is inspired by the fact that $g^\flat$ lowers indices of tensors, wheras $g^\sharp$ raises them.

\begin{definition}{Remark}
As it is well known, if $(M,g)$ is a Riemannian supermanifold of graded dimension $p|q$, then the odd dimension $q$ must be even because of the non-degeneracy condition of the super metric.
Note that the definition of a super metric as given here is the straightforward generalisation of a metric to the super context. In \cite[\S IV.7]{Tu1} a different (and not completely natural) notion of a super metric was introduced. That definition was adapted to the need to be able to define a supplement to any subbundle of a given vector bundle without the constraint that the odd dimension should be even. 
%\textcolor{blue}{J'insiste, car en tant qu'auteur de mon livre je ne peux pas ne pas faire cette remarque. Sinon on pourrait croire que j'ai changé d'avis sur la définition de super métrique.}%

\end{definition}

%If $(x^1, \dots, x^n)$ are local coordinates on $M$, remember that the vectors $\partial_{x^i}\vert_m$ form a basis of the tangent space to $T_mM$. We define a matrix $g_{ij}$ by

If $(x^1, \dots, x^n)$ are local coordinates on $M$, then the vectors $\partial_{x^i}\vert_m$ form a basis of the tangent space $T_mM$. %\textcolor{blue}{au choix : ``the tangent space $T_mM$'' ou ``the tangent space to $M$ at $m$'', mais jamais ``the tangent space to $T_mM$''.}%
 Using these vectors, we define the matrix $g_{ij}$ by

%\begin{moneq}
%g_{ij} = g_m(\partial_{x^i}\vert_m,\partial_{x^j}\vert_m)
%\end{moneq}%
%\textcolor{green}{AU CHOIX}
\begin{moneq}
g_{ij} = \contrf{\partial_{x^i}\vert_m,\partial_{x^j}\vert_m}{g_m}
\end{moneq}%
It follows immediately that for any two arbitrary tangent vectors $v = \sum_i v^i \,\partial_{x^i}\vert_m$ and $w = \sum_i w^i \,\partial_{x^i}\vert_m$, we have
%\begin{moneq}
%g_m(v,w)
%=
%\sum_{i,j} v^i \, \conjug^{\varepsilon_i}(w^j) \, g_{ij}
%\end{moneq}
%\textcolor{green}{AU CHOIX}
\begin{moneq}
\contrf{v,w}{g_m}
=
\sum_{i,j} v^i \, \conjug^{\varepsilon_i}(w^j) \, g_{ij}
\end{moneq}%
Equivalently, in terms of the (left-)dual basis $(\extder x^1\vert_m, \dots, \extder x^n\vert_m)$ of $T_m^*M$, we have %\textcolor{blue}{si on écrit ``we write'', pour moi ça veut dire que quelqu'un d'autre peut écrire autre chose, ce qui n'est pas le cas.}%
\begin{moneq}
g_m
=
\sum_{ij} \extder x^j\vert_m \otimes \extder x^i\vert_m \, g_{ij}
\end{moneq}%
The graded-symmetry and even-ness of $g_m$ translate as %\textcolor{blue}{vous avez remplacé mon ``as'' par ``into'', ce qui est faux : on dit ``the word building is translated \textcolor{green}{into} french \textcolor{green}{as} bâtiment''}%
 the properties
\begin{moneq}%[gradedsymofg]
g_{ij} = (-1)^{\varepsilon_i\,\varepsilon_j} \, g_{ji}
\qquad\text{and}\qquad
\varepsilon(g_{ij}) = \varepsilon_i + \varepsilon_j
\end{moneq}%
and
%Also, the 
non-degeneracy means that the matrix $g_{ij}$ is invertible. We denote the inverse matrix by $g^{ij}$, \ie, we have the equalities
\begin{moneq}%[defofraisingindicesop]
\sum_j g_{ij} \,g^{jk} = \delta_i^k = \sum_j g^{kj} \, g_{ji}
\end{moneq}%
where $\delta_i^k$ denotes the Kronecker delta. It is straightforward that the parity of $g^{ij}$ is $\varepsilon(g^{ij}) = \varepsilon_i + \varepsilon_j$, while the graded symmetry of $g$ gives us the following symmetry property of the inverse matrix:
\begin{moneq}%[antigradedsymofgmo]
g^{ij} = (-1)^{\varepsilon_i + \varepsilon_j + \varepsilon_i\varepsilon_j} g^{ji}
\end{moneq}%
Finally note that the map $g_m^\flat:T_mM \to T_m^*M$ reads
%\begin{moneq}
%g_m^\flat(v)
%=
%\sum_{ij} (-1)^{\varepsilon_i} \, v^j \, g_{ji} \, \extder x^i\vert_m
%\qquad\text{ for $v = \sum_i v^i \,\partial_{x^i}\vert_m$}
%\end{moneq}
%\textcolor{green}{AU CHOIX}
\begin{moneq}
\contrf{v}{g_m^\flat}
=
\sum_{ij} (-1)^{\varepsilon_i} \, v^j \, g_{ji} \, \extder x^i\vert_m
\qquad\text{ for $v = \sum_i v^i \,\partial_{x^i}\vert_m$}
\end{moneq}
and that, using the inverse matrix, it is not hard to show that 
%\textcolor{blue}{j'ai enlevé du text}%
%the value on an arbitrary $\alpha = \sum_i \alpha_i \extder x^i\vert_m \in T_m^*M$ of 
the inverse map $g_m^\sharp = (g_m^\flat)^{-1} : T_m^*M \to T_mM$ is given by
%\begin{align}
%g_m^\sharp(\alpha)
%=
%\sum_{ij}
%(-1)^{\varepsilon_i} \,\alpha_i \, g^{ij} \, \partial_{x^j}
%\qquad\text{ for $\alpha = \sum_i \alpha_i \,\extder x^i\vert_m$}
%\label{gsharpinlocalcoord}
%\end{align}
%\textcolor{green}{AU CHOIX}
\begin{align}
\contrf{\alpha}{g_m^\sharp}
=
\sum_{ij}
(-1)^{\varepsilon_i} \,\alpha_i \, g^{ij} \, \partial_{x^j}
\qquad\text{ for $\alpha = \sum_i \alpha_i \,\extder x^i\vert_m$}
\label{gsharpinlocalcoord}
\end{align}%

\begin{proclaim}[EandUofmetricconnection]{Lemma}
If $(M,g)$ is a Riemannian supermanifold, there exists a unique torsion-free connection $\nabla$ in $TM$ which is compatible with the metric in the sense that for any three homogeneous vector fields $X$, $Y$ and $Z$ on $M$, we have
%\begin{moneq}[compatibleconnection]
%Xg(Y,Z)
%=
%g(\nabla_X Y, Z) + (-1)^{\varepsilon(X)\varepsilon{(Y)}} \, g(Y, \nabla_XZ)
%\end{moneq}
%\textcolor{green}{AU CHOIX}
\begin{moneq}[compatibleconnection]
X \bigl(\, \contrf{Y,Z}{g} \,\bigr)
=
\contrf{\nabla_X Y, Z}{g} + (-1)^{\varepsilon(X)\varepsilon{(Y)}} \, \contrf{Y, \nabla_XZ}{g}
\end{moneq}

\end{proclaim}

\begin{preuve}
Existence follows from the explicit formula for the Christoffel symbols in local coordinates
\begin{moneq}
\Gamma_{jk}^i
=
\tfrac12 \sum_\ell\bigl(\ 
\partial_{x^j}g_{k\ell} + (-1)^{\varepsilon_{j} \varepsilon_{k}}\, \partial_{x^k} g_{j\ell} - (-1)^{\varepsilon_{\ell}(\varepsilon_{j} + \varepsilon_{k})}\,\partial_{\ell} g_{jk}
\bigr)\, g^{\ell i}
\end{moneq}%
For uniqueness we observe first that condition \recalle{compatibleconnection} applied to the (local) vector fields $X=\partial_{x^p}$, $Y=\partial_{x^j}$ and $Z=\partial_{x^k}$ gives us the equality
\begin{moneq}
\partial_{x^p} g_{jk}
=
\Gamma_p{}^i{}_j \,g_{ik} + (-1)^{\varepsilon_j\varepsilon_k}\, \Gamma_p{}^i{}_k \, g_{ij}
\end{moneq}%
It follows that if we have two connections $\nabla$ and $\nablah$ satisfying these conditions, then the components $S_{jk}^i = \Gamma_{jk}^i - \Gammah_{jk}^i$ of the difference tensor must satisfy the conditions
\begin{moneq}
S_{pj}^i \,g_{ik} = -(-1)^{\varepsilon_j\varepsilon_k} \,S_{pk}^i \,g_{ij}
\end{moneq}%
Using the graded symmetry of the tensor $S$ (the connections are torsion-free), we can further compute
\begin{align*}
S_{pj}^i\, g_{ik}
&=
(-1)^{\varepsilon_j\varepsilon_p}\, S_{jp}^i \, g_{ik}
=
(-1)^{1+\varepsilon_p(\varepsilon_j + \varepsilon_k)} \,S_{jk}^i g_{ip}
\\&
=
(-1)^{1+\varepsilon_p(\varepsilon_j + \varepsilon_k) + \varepsilon_j\varepsilon_k} \,S_{kj}^i g_{ip}
=
(-1)^{\varepsilon_k(\varepsilon_j + \varepsilon_p)} \,S_{kp}^i g_{ij}
\\&
=
(-1)^{\varepsilon_k\varepsilon_j } \,S_{pk}^i g_{ij}
=
- S_{pj}^i g_{ik}
\end{align*}
This shows that the difference tensor must be zero, \ie, $\nabla=\nablah$.
\end{preuve}

\begin{definition}{Definition}
Let $\pr:T^*M\to M$ be the cotangent bundle of the supermanifold $M$. The \emph{canonical $1$-form} $\theta$ on $T^*M$ is defined as follows: for $\alpha\in T^*M$ and $V \in T_\alpha(T^*M)$ we write $m=\pr(\alpha)$ (and thus $\alpha\in T_m^*M$), and then
\begin{moneq}
\contrf{V}{\theta_\alpha} = \contrf{v}{\alpha}
\end{moneq}%
where $v = \contrf{V}{T\pr} \in T_mM$ %\textcolor{blue}{j'ai changé l'ordre, car les vecteurs tangent sont à gauche de l'application}%
 is the image of $V \in T_\alpha(T^*M)$ under the tangent map of the canonical projection.
%by the tangent map of $\pr$.
\end{definition}

If if $(x^1, \dots, x^n)$ are local coordinates on $M$, then any $1$-form $\alpha$ at $m\in M$ can be expressed as $\alpha = \sum_i \alpha_i \, \extder x^i$. Splitting the coefficients $\alpha_i \in \CA$ into their even and odd parts $\alpha_i = p_i + \pb_i$, we write
\begin{moneq}
\alpha = \sum_i (p_i + \pb_i) \, \extder x^i
\qquad\text{with}\qquad
\alpha_0 = \sum_i p_i \, \extder x^i
\quad\text{and}\quad
\alpha_1 = \sum_i \pb_i \,\extder x^i
\end{moneq}%
The parity of these coordinates thus is given by $\varepsilon(p_i) = \varepsilon_i$ and $\varepsilon(\pb_i) = \varepsilon_i +1$. Thus, if the graded dimension of $M$ is $p\vert q$, then the graded dimension of the full cotangent bundle is $2p+q\vert p+2q$ with coordinates $x^i$, $p_i$ and $\pb_i$, the graded dimension of its even part (whose sections are the even $1$-forms) is $2p\vert2q$ with coordinates $x^i$ and $p_i$ and the graded dimension of its odd part (whose sections are the odd $1$-forms) is $p+q\vert p+q$ with coordinates $x^i$ and $\pb_i$.

In terms of these local coordinates on $T^*M$, it is easy to show that the canonical $1$-form $\theta$ on $T^*M$ is given by
\begin{moneq}
\theta = \sum_i (p_i+\pb_i) \,\extder x^i
\end{moneq}%
By definition, the \emph{canonical $2$-form} $\omega$ on $T^*M$ is the exterior derivative of the canonical $1$-form: $\omega = \extder \theta$. In local coordinates $\omega$ thus reads
\begin{moneq}
\omega = \sum_i \extder p_i \wedge \extder x^i + \sum_i \extder \pb_i \wedge \extder x^i
\end{moneq}%
In particular, the restriction of $\omega$ to $T^*M^{(0)}$, the even part of the cotangent bundle, is an even symplectic form, while its restriction to the odd part of the cotangent bundle $T^*M^{(1)}$ is an odd symplectic form. 

%\textcolor{blue}{j'ai changé un tout petit peu, car il n'y a pas trois cotangent différents ($T^*M$, le cotangent pair et le cotangent impair), mais seulement un, dont on prend deux ``sous-fibrés'' (pas super vectoriel !!)}%

\medskip

%\begin{definition}{Remark}
%Following \cite{Tu2} we observe that the non-homogeneous $2$-form $\omega$ on the full cotangent bundle is degenerate, but homogeneously non-degenerate. Hence the full cotangent bundle is a symplectic supermanifold with a non-homogeneous symplectic form. However, as symplectic manifold it is not very interesting because the Poisson algebra $\mathcal{P}(T^*M) $ turns out to be rather ``small''. Roughly speaking, $\mathcal{P}(T^*M)$ consists in pairs of functions that are affine in the cotangent coordinates $p_i, \pb_i$. More precisely, one finds
%\begin{multline*}
%\mathcal{P}(T^*M) 
%=
%\{\, (f_0,f_1)\in C^\infty(T^*M)^2 \mid
%\exists f_{00}, f_{10} \in C^\infty(M) \ ,\ \exists X\in \Gamma(TM) : 
%\\
%f_0(\alpha) = f_{00}(\pr(\alpha)) + \contrf{X}{\alpha_0}
%\ ,\ f_1(\alpha) = f_{10}(\pr(\alpha) + \contrf{X}{\alpha_1}) 
%\,\}
%\end{multline*}%
%where $\Gamma(TM)$ stands for the $\RR$-module of smooth vector fields on $M$.
%\end{definition}

We now come to the description of the movement of a free particle with unit mass on the Riemannian supermanifold $(M,g)$. There is no potential energy while kinetic energy is simply given by half the metric. More precisely, the phase space is the even part of the cotangent bundle $T^*M^{(0)}$ while the \emph{Hamiltonian} of the system is the function $H:T^*M^{(0)}\to \CA$ whose value on an element $\alpha\in T_m^*M^{(0)}$ is 
\begin{moneq}[kinetichamiltonian]
H(\alpha) = \tfrac12 \, \contrf{g_m^\sharp(\alpha), g_m^\sharp(\alpha)}{g_m}
%\equiv \tfrac12 \, g_m(\,g_m^\sharp(\alpha), g_m^\sharp(\alpha) \, )
\end{moneq}%
In local coordinates, the Hamiltonian thus reads
\begin{align*}
H(x,p) 
&= 
\tfrac12\sum_{jk} (-1)^{\varepsilon_j+ \varepsilon_k} \, p_j \, g^{jk}(x)  \, p_k
=
\tfrac12\sum_{jk} (-1)^{\varepsilon_j} \, p_k \, p_j \, g^{jk}(x) 
\\&
=
\tfrac12\sum_{jk} (-1)^{\varepsilon_k} \, g^{jk}(x) \, p_k \, p_j  
\end{align*}%
The local expression for $\omega$ is $\omega = \sum_i \extder p_i \wedge \extder x^i$ and the definition of the hamiltonian vector field $X_f$ associated with a function $f$ is given by the formula
\begin{moneq}%[defhamvectfield]
\contrf{X_f}{\omega} = - \extder f
\end{moneq}%
In local coordinates this gives us
\begin{moneq}
X_f
=
\sum_i 
\bigl(\,
(-1)^{\varepsilon_i} \, \conjug^{\varepsilon_i}(\partial_{p_i}f) \,\partial_{x^i} - \conjug^{\varepsilon_i}(\partial_{x^i}f) \, \partial_{p_i}
\,\bigr)
\end{moneq}%
and thus, for our particular function $H$, we obtain the even vector field
\begin{moneq}
X_H
=
\sum_{ik} (-1)^{\varepsilon_k} \,p_k \, g^{ki} \, \fracp{}{x^i}
-
\tfrac12\,\sum_{ijk} (-1)^{\varepsilon_i +  \varepsilon_k } \, \fracp{g^{jk}}{x^i} \, p_k \, p_j \, \fracp{}{p_i}
\end{moneq}

%\begin{definition}{Remark}
%Knowing that we also have a symplectic form on the odd tangent bundle and on the full tangent bundle, we could have tried to play the same game on these symplectic manifolds. However, formula \recalle{kinetichamiltonian} applied to elements of $T^*M^{(1)}$ gives us a function which is identically zero, simply because $g$ is graded symmetric and $g_m^\sharp(\alpha)$ is an odd tangent vector. So on the odd tangent bundle nothing interesting happens. Also, formula \recalle{kinetichamiltonian} applied to the full tangent bundle yields a function which is quadratic in the momentum coordinates, and thus no splitting into two functions can give us an element of the Poisson algebra. So again nothing interesting is obtained.
%\end{definition}

\begin{definition}{Remark}
Knowing that we also have a symplectic form on the odd tangent bundle and on the full tangent bundle, we could have tried to play the same game on these symplectic manifolds. However, formula \recalle{kinetichamiltonian} applied to elements of $T^*M^{(1)}$ gives us a function which is identically zero, simply because $g$ is graded symmetric and $g_m^\sharp(\alpha)$ is an odd tangent vector. So on the odd tangent bundle nothing interesting happens. 
Note that the full cotangent bundle is also a symplectic supermanifold (with a non-homogeneous symplectic form). However, it can be shown following \cite{Tu2} that formula \recalle{kinetichamiltonian} yields a function which is not in the Poisson algebra of $T^*M$, \ie, a function which does not give rise to a hamiltonian vector field. So again nothing interesting can be obtained.
\end{definition}

\begin{proclaim}{Proposition}
Under the isomorphism $g^\sharp : T^*M^{(0)} \to TM^{(0)}$ the vector field $X_H$ on $T^*M^{(0)}$ is mapped to the vector field $G$ on $TM^{(0)}$ given by \recalle{definitionofG} using the unique metric connection given by \recallp{EandUofmetricconnection}
\end{proclaim}

\begin{preuve}
The proof is a lenghty but straightforward computation.
\end{preuve}

It follows that the integral curves of the Hamiltonian vector field $X_H$ correspond to the integral curves of the geodesic vector field of the metric connection associated with $g$, and thus in particular the geodesics of the metric connection coincide with the projections of the integral curves of the Hamiltonian vector field onto $M$, \ie, the geodesics are the trajectories of a free particle with unit mass on the Riemannian supermanifold $(M,g)$.

\begin{definition}{Remarks}

$\bullet$
The isomorphism $g^\sharp : T^*M^{(0)} \to TM^{(0)}$ can be interpreted as the Legendre transformation, which transforms the Hamiltonian formalism on the cotangent bundle into the Lagrangean formalism on the tangent bundle. More details on this interpretation in the non-super case can be found in \cite[\S3.6--7]{AM}.

$\bullet$
We have used {\em left} coordinates $p_i, \pb_i$ on the cotangent bundle, writing $\alpha = \sum_i (p_i+\pb_i) \,\extder x^i$. We could also have used {\em right} coordinates $p_i', \pb'_i$ by writing $\alpha = \sum_i \extder x^i \, (p'_i + \pb'_i)$. They are related by the simple equations $\pb'_i = \pb_i$ and $p'_i = (-1)^{\varepsilon_i} \, p_i$. 
This would have ``simplified'' the formul{\ae} for $H$ to 

\begin{moneq}
H(x,p') = \tfrac12\, \sum_{jk}p'_j\, g^{jk} \, p'_k
\end{moneq}%
The reason not to use these coordinates (and it {\em is} a simple change of coordinates) is first that it is good practice not to mix left- and right-coordinates at the same time (and when using matrices it becomes crucial, see \cite[VI.1.20]{Tu1}) and secondly that the explicit expression for the full map $g^\sharp:T^*M\to TM$ would have contained the conjugation map $\conjug$, as we would have had to transform the right-coordinates $\alpha_i$ of $\alpha=\sum_i \extder x^i \, \alpha_i$ into left coordinates $v^j$ of $v=\sum_j v^j \, \partial_{x^j} = g^{\sharp}(\alpha)$.
\end{definition}

\appendix

\section{The exponential map}

In the non-super case it is well known that ``running faster'' through a geodesic is the same as taking the geodesic with a bigger initial velocity. In terms of the flow $\Psi\cong(\Psi_1,\Psi_2)$ this would mean that we should have
\begin{moneq}
\Psi_1(t,x,\lambda v) = \Psi_1(\lambda t,x,v)
\qquad\text{and}\qquad
\Psi_2(t,x,\lambda v) = \lambda\cdot \Psi_2(\lambda t, x,v)
\end{moneq}%
for any $\lambda\in \CA_0$. 

In order to prove this rigourously and in a coordinate independent way, we introduce the map $D_\lambda:TM^{(0)}\to TM^{(0)}$, the dilation of the tangent space by a factor $\lambda$, in local coordinates by
\begin{moneq}
D_\lambda(x,v) = (x,\lambda v)
\end{moneq}%
These local definitions glue together to form a well-defined global map. Moreover, it does not affect the base point:
\begin{moneq}
\pi \scirc D_\lambda = \pi : TM^{(0)}\to M
\end{moneq}%

\begin{proclaim}{Proposition}
On a suitable open domain in $\CA_0\times\CA_0\times TM^{(0)}$ containing $\{0\}\times \{0\}\times TM^{(0)}$, the maps $\Psih$ and $\Psit$ with values in $TM^{(0)}$ and defined by
\begin{align*}
\Psit(t,\lambda,\tanvec) 
&=
\Psi(t, D_\lambda(\tanvec)) \cong (\Psi_1(t,x,\lambda v), \Psi_2(t,x,\lambda v))
\\
\Psih(t,\lambda,\tanvec) 
&=
D_\lambda( \Psi(\lambda t, \tanvec)) \cong (\Psi_1(\lambda t, x, v), \lambda \cdot \Psi_2(\lambda t,x,v))
\end{align*}
are the same.
\end{proclaim}

\begin{preuve}
We start with the observation that in local coordinates $(x,v)$ on $TM^{(0)}$ the tangent map of $D_\lambda$ behaves as
\begin{moneq}
\contrf{\partial_{x^i}\vert_{(x,v)}}{TD_\lambda}
=
\partial_{x^i}\vert_{(x,\lambda v)}
\qquad\text{and}\qquad
\contrf{\partial_{v^i}\vert_{(x,v)}}{TD_\lambda}
=
\lambda\cdot \partial_{v^i}\vert_{(x,\lambda v)}
\end{moneq}%
It follows that we have the following equality concerning the local expression of the vector field $G$:
\begin{align*}
\lambda\cdot \contrf{G\vert_{(x,v)}}{TD_\lambda}
&=
\lambda\cdot\contrf{\sum_i v^i \partial_{x^i}\vert_{(x,v)}
-
\sum_{ijk}
v^k\cdot v^j\cdot\Gamma_{jk}^i(x) 
\cdot \partial_{v^i}\vert_{(x,v)}}{TD_\lambda}
\\&
=
\lambda\cdot\sum_i v^i \partial_{x^i}\vert_{(x,\lambda v)}
-+
\sum_{ijk}
\lambda\cdot v^k\cdot v^j\cdot\Gamma_{jk}^i(x) 
\cdot\lambda\cdot \partial_{v^i}\vert_{(x,\lambda v)}
\\&
=
G\vert_{(x,\lambda v)}
\end{align*}%
which means that $\lambda\cdot G\vert_{\tanvec}$ is mapped by $TD_\lambda$ to $G\vert_{D_\lambda(\tanvec)}$.

With that knowledge we compute the image of the tangent vector $\partial_t$ under the maps $\Psit$ and $\Psih$:
\begin{subequations}
\label{PsihendPsitimageofpartialt}
\begin{moneq}[aap]
\contrf{\partial_t\vert_{(t,\lambda,\tanvec)}}{T\Psit}
=
\contrf{\partial_t\vert_{(t,D_\lambda(\tanvec))}}{T\Psi}
=
G\vert_{\Psi(t,D_\lambda(\tanvec))}
= 
G\vert_{\Psit(t,\lambda,\tanvec)}
\end{moneq}%
and
\begin{align}
\notag
\contrf{\partial_t\vert_{(t,\lambda,\tanvec)}}{T\Psih}
&=
\lambda\cdot\contrf{\partial_t\vert{(\lambda t,\tanvec)}}{T(D_\lambda\scirc \Psi)}
=
\lambda\cdot\contrf{G\vert_{\Psi(\lambda t, \tanvec)}}{TD_\lambda}
\\&
=
G\vert_{D_\lambda(\Psi(\lambda t, \tanvec))}
=
G\vert_{\Psih(t,\lambda, \tanvec)}
\end{align}
\end{subequations}

We then introduce the extended manifold $N= \CA_0 \times TM^{(0)}$ on which we define the even vector field $H$ (the extension of $G$ to $N$) by
\begin{moneq}
H\vert_{(\lambda,\tanvec)}
=
G\vert_{\tanvec}
\end{moneq}%
and we introduce the maps $\Phit, \Phih:\CA_0 \times N \to N$ by
\begin{moneq}
\Phit(t,\lambda,\tanvec)
=
\bigl(\lambda, \Psit(t,\lambda, \tanvec)\bigr)
\qquad\text{and}\qquad
\Phih(t,\lambda,\tanvec)
=
\bigl(\lambda, \Psih(t,\lambda, \tanvec)\bigr)
\end{moneq}%
It then is immediate from \recalle{PsihendPsitimageofpartialt} that we have
\begin{moneq}
\contrf{\partial_t\vert_{(t,\lambda,\tanvec)}}{T\Phit}
=
H\vert_{\Phit(t,\lambda,\tanvec)}
\qquad\text{and}\qquad
\contrf{\partial_t\vert_{(t,\lambda,\tanvec)}}{T\Phih}
=
H\vert_{\Phih(t,\lambda,\tanvec)}
\end{moneq}%
Moreover, at time $t=0$ we have
\begin{moneq}
\Phit(0,\lambda,\tanvec) 
= 
\bigl(\lambda,D_\lambda(\tanvec)\bigr)
=
\Phih(0,\lambda,\tanvec)
\end{moneq}%
As the map $(\lambda,\tanvec) \mapsto D_\lambda(\tanvec)$ is smooth, we can apply the (existence and) uniqueness of local flows of a vector field ($H$ in our case) with given initial condition to conclude that $\Phit$ and $\Phih$ and thus a fortiori $\Psit$ and $\Psih$ are the same \cite[V.4.8]{Tu1}.
\end{preuve}

\begin{definition}{Remark}
We have been a bit vague on the domain of definition on which the maps are defined. The domains of $\Psit$ and $\Psih$ are in the obvious way related to the domain $W_G$ of the flow $\Psi$, but initially it is not clear that they are the same. The fact that these two maps coïncide then proves that these two domains coïncide. And thus that we have in particular the equivalence
\begin{moneq}
(\lambda t, \tanvec) \in W_G
\quad\Longleftrightarrow\quad
\bigl(t,D_\lambda(\tanvec)\bigr) \in W_G
\end{moneq}%

\end{definition}

\begin{proclaim}{Corollary}
Running faster through a geodesic is the same as taking a bigger initial velocity:
\begin{moneq}
\pi( \Psi(\lambda t,\tanvec))
=
\pi( \Psi(t,D_\lambda(\tanvec)))
\end{moneq}%
In local coordinates this boils down to $\Psi_1(\lambda t,x,v) = \Psi_1(t,x,\lambda v)$. Moreover, the subset $\Omega\subset TM^{(0)}$ defined as
\begin{moneq}
\Omega = \{\,\tanvec \in TM^{(0)} \mid (1,\tanvec) \in W_G \,\}
\end{moneq}%
contains the zero section of the tangent bundle $TM^{(0)}$.

\end{proclaim}

\begin{definition}{Definition}
Let $\nabla$ be a connection on $TM$ and let $\Psi:W_G\to TM^{(0)}$ be the flow of the vector field $G$ associated with $\nabla$. Then the geodesic exponential map $\exp:\Omega \to M$ is defined as
\begin{moneq}
\tanvec\in T_mM^{(0)} \mapsto \exp_m(\tanvec) =
\pi\bigl( \Psi(1,\tanvec)\bigr)
\qquad\text{with $m=\pi(\tanvec)$}
\end{moneq}%
This map is jointly smooth in the coordinates $(x,v)$ of $\tanvec\in \Omega$. However, if $m=\pi(\tanvec)$ does not belong to the body of $M$, then there is no guarantee that the map $\exp_m : T_mM^{(0)} \to M$ (with $m$ fixed) is smooth.
\end{definition}

% SECTION
\section*{Acknowledgments}

It is a pleasure to thank S.~Garnier for stimulating discussions. F. Radoux also thanks the Belgian FNRS for his research fellowship.

\end{document}